\documentclass[a4paper,12pt,numbers,sort&compress]{article}

\pdfoutput=1

\usepackage[paper=letterpaper,margin=1in]{geometry}

\usepackage{amsmath,amssymb,amsthm,amsfonts,epsfig,cite,setspace,bigstrut,longtable,array,breqn,url,color}

\usepackage{caption,setspace}

\usepackage{tikz}
\usepackage{mathrsfs}
\usetikzlibrary{arrows}

\usepackage{pdfpages,lipsum}

\begin{document}



\vskip 0.25in

\newcommand{\todo}[1]{{\bf ?????!!!! #1 ?????!!!!}\marginpar{$\Longleftarrow$}}
\newcommand{\sref}[1]{\S~\ref{#1}}
\newcommand{\nn}{\nonumber}
\newcommand{\tr}{\mathop{\rm Tr}}
\newcommand{\comment}[1]{}

\newcommand{\cM}{{\cal M}}
\newcommand{\cW}{{\cal W}}
\newcommand{\cN}{{\cal N}}
\newcommand{\cH}{{\cal H}}
\newcommand{\cK}{{\cal K}}
\newcommand{\cZ}{{\cal Z}}
\newcommand{\cO}{{\cal O}}
\newcommand{\cP}{{\cal P}}
\newcommand{\cR}{{\cal R}}
\newcommand{\cA}{{\cal A}}
\newcommand{\cB}{{\cal B}}
\newcommand{\cC}{{\cal C}}
\newcommand{\cD}{{\cal D}}
\newcommand{\cE}{{\cal E}}
\newcommand{\cF}{{\cal F}}
\newcommand{\cT}{{\cal T}}
\newcommand{\cV}{{\cal V}}
\newcommand{\cX}{{\cal X}}
\newcommand{\IA}{\mathbb{A}}
\newcommand{\IP}{\mathbb{P}}
\newcommand{\IQ}{\mathbb{Q}}
\newcommand{\IH}{\mathbb{H}}
\newcommand{\IK}{\mathbb{K}}
\newcommand{\IR}{\mathbb{R}}
\newcommand{\IC}{\mathbb{C}}
\newcommand{\IF}{\mathbb{F}}
\newcommand{\IV}{\mathbb{V}}
\newcommand{\II}{\mathbb{I}}
\newcommand{\IZ}{\mathbb{Z}}
\newcommand{\re}{{\rm~Re}}
\newcommand{\im}{{\rm~Im}}

\newcommand{\tmat}[1]{{\tiny \left(\begin{matrix} #1 \end{matrix}\right)}}
\newcommand{\mat}[1]{\left(\begin{matrix} #1 \end{matrix}\right)}

\let\oldthebibliography=\thebibliography
\let\endoldthebibliography=\endthebibliography
\renewenvironment{thebibliography}[1]{%
\begin{oldthebibliography}{#1}%
\setlength{\parskip}{0ex}%
\setlength{\itemsep}{0ex}%
}%
{%
\end{oldthebibliography}%
}

\newtheorem{theorem}{\bf THEOREM}
\newtheorem{proposition}{\bf PROPOSITION}
\newtheorem{observation}{\bf OBSERVATION}
\newtheorem{definition}{\bf DEFINITION}
\def\theequation{\thesection.\arabic{equation}}

\newcommand{\setall}{\setcounter{equation}{0}}

\begin{titlepage}

~\\
\vskip 1cm
\begin{center}
{\Large \bf Graph Laplacians, Riemannian Manifolds \\
\qquad \& their Machine-Learning}

\medskip

Yang-Hui He$^{1,2,3}$ \&
Shing-Tung Yau$^{4,5,6}$

\renewcommand{\arraystretch}{0.5} 
{\small
{\it
\begin{tabular}{rl}
  ${}^{1}$ &
  Merton College, University of Oxford, OX14JD, UK\\
  ${}^{2}$ &
    Department of Mathematics, City, University of London, EC1V 0HB, UK\\
  ${}^{3}$ &
    School of Physics, NanKai University, Tianjin, 300071, P.R.~China\\
  ${}^{4}$ &
  Department of Mathematics,  Department of Physics, \& Center of Mathematical Sciences \\
  & ~~ and Applications, Harvard University, Cambridge, MA 02138, USA\\
  ${}^{5}$ &
  Yau Mathematical Sciences Center, Tsinghua University, Beijing, 100804, China\\
  ${}^{6}$ &
  Beijing Institute of Mathematical Sciences and Applications, \\
  & ~~ Huairou Science City, Beijing, 101400, China\\
\end{tabular}
}
~\\
~\\
~\\
hey@maths.ox.ac.uk; \ yau@math.harvard.edu
}
\renewcommand{\arraystretch}{1.5} 

\end{center}

\vspace{10mm}

\begin{abstract}
Graph Laplacians as well as related spectral inequalities and (co-)homology provide a foray into discrete analogues of Riemannian manifolds, providing a rich interplay between combinatorics, geometry and theoretical physics.
We apply some of the latest techniques in data science such as supervised and unsupervised machine-learning and topological data analysis to the Wolfram database of some 8000 finite graphs in light of studying these correspondences.
Encouragingly, we find that neural classifiers, regressors and networks can perform, with high efficiently and accuracy,
a multitude of tasks ranging from recognizing graph Ricci-flatness, to predicting the spectral gap, to detecting the presence of Hamiltonian cycles, etc.
\end{abstract}

\end{titlepage}

\tableofcontents

\section{Introduction and Summary}\setall
The Laplacian is of paramount importance in mathematics, its ubiquity has extended from curvature in differential geometry, to (co-)homology in algebraic geometry, to invariants in topology and to particle spectrum in high energy physics.
In classical Riemannian geometry, the celebrated decomposition of Hodge relates the zero-modes of the Laplacian to the representative in cohomology, and hence allows for an elegant computation $H^*$.
In complex geometry, the second author's reduction of the Laplacian eigen-equation for K\"ahler manifolds to analyses of a PDE of Monge-Amp\`ere type \cite{yau} resulted in the  proof of Calabi's Conjecture \cite{calabi}.

A natural question arose as to whether there should exist a {\it discrete} version of the story.
A programme had been launched over the last decade or so to understand it in the context of locally finite graphs.
Based on the standard theory of the Laplace operator on graphs (cf.~ e.g., \cite{FanBook}), the authors in \cite{fy,lly,ly2,HuaLin,glly,ckllly} pursued a notion of curvature on graphs to examine the analogue for manifolds
(q.v.~a comprehensive review in \cite{ly}). 
Furthermore, \cite{glmy1,glmy2} investigated the idea of (co-)homology of finite graphs.

In parallel, a recent programme of using latest techniques from machine-learning and data science to study various mathematical formulae and conjectures  had been proposed \cite{He:2017aed,He:2018jtw}.
Experimentation of whether standard techniques in neural networks and classifiers could be carried over to study diverse problems in mathematics have ranged from triangulations in Calabi-Yau hypersurfaces in toric varieties \cite{Altman:2018zlc,Demirtas:2018akl,He:2015fif}, to finding bundle cohomology on varieties \cite{Ruehle:2017mzq,Brodie:2019dfx,Larfors:2020ugo}, to distinguishing elliptic fibrations \cite{Anderson:2017aux,He:2019vsj} and invariants of Calabi-Yau threefolds \cite{Bull:2018uow} (cf.~\cite{Grimm:2019bey,Grimminger:2020dmg} on the organization of classification by Hasse diagrams), to knot hyperbolic volumes \cite{Jejjala:2019kio}, to machine-learning the Donaldson algorithm for numerical Calabi-Yau metrics \cite{Ashmore:2019wzb}, to the algebraic structures of groups and rings \cite{He:2019nzx}, to the BSD conjecture in number theory \cite{Alessandretti:2019jbs}, to finding discriminant locii \cite{discriminant} etc. (q.v.~\cite{hetalk} for speculations on how the foundations of mathematics might respond to machine-learning).

Indeed, with the introduction of the machine-learning paradigm \cite{He:2017aed,Krefl:2017yox,Ruehle:2017mzq,Carifio:2017bov} to string theory, a multitude of heartening results in physics have included, e.g., finding Higgsable gauge groups \cite{Wang:2018rkk}, axion physics from strings \cite{Demirtas:2018akl}, flux compactifications \cite{Cole:2019enn}, distinguishing standard models \cite{Mutter:2018sra,Otsuka:2020nsk,Deen:2020dlf}, QFT dualities \cite{Betzler:2020rfg,Bao:2020nbi,He:2020eva}, detecting symmetries \cite{Krippendorf:2020gny} and CFTs \cite{1803639}, etc.
Of particular note is  \cite{Hashimoto:2018ftp} where the holographic AdS/CFT correspondence and thence, space-time  itself, is interpreted as a neural-network/Boltzmann machine.
The reader is furthermore referred to fascinating works in the last couple of years on neural-networks which can perform symbolic mathematics \cite{lc}, find results in fundamental physics \cite{imwdr,Cranmer:2020wew} from scratch and discover chemical reactions \cite{chem} from word-embedding (cf. a similar linguistic study \cite{He:2018dlv} of ArXiv titles in classifying different disciplines in mathematical physics and in generating syntactical identities).

Given the highly structured representation of finite graphs in terms of matrix manipulations - something for which machine-learning is perfectly adapted - it is immediate to ask whether the aforementioned two programmes should intersect.
That is, whether machine-learning could be applied to the investigation of the geometry of graphs.
It is this question which we will address in this paper, which will hopefully initiate a novel direction in the study of the geometry and algebra in graph theory, especially in the context of discrete analogues to Riemannian manifolds.

Bearing this question in mind, the paper is organized as follows.
In Section 2, we begin with some rudiments on the graph Laplacian and related standard facts in graph theory, as well as the connection to Riemannian geometry.
Then, in Section 3, we present our main protagonist of the Wolfram Database of finite, undirected, simple graphs, a set of some 8000 graphs organized by vertices and edges, and named wherever possible.
On this set we examine some preliminary properties, such as planarity, genus, chromatic number, etc., to see whether they can be machine-learned.

Section 4 is then concerned with spectral bounds. These are bounds on the eigenvalues of the Laplacian, which are not only important to the classical spectral theory of graphs, but are also enlightening to Li-Yau type of inequalities \cite{LiYau} in compact manifolds. We study the distribution of the eigenvalues and then apply topological data analysis, principal component analysis, as well as supervised machine-learning to these bounds. 
Section 5 is devoted to the critical class of finite graphs which are Ricci-flat.
These graphs are the analogues of Calabi-Yau manifolds, and the classification thereof has been an active field of research \cite{lly,ly2,HuaLin}. We will show that a neural classifier can distinguish a Ricci-flat graph to very high accuracy.

Finally, Section 6 is devoted to studying the homology of directed (finite, simple) graphs.
We explicitly compute the Euler number of the graphs in the sense of \cite{glmy1,glmy2} for the Wolfram database and see how machine-learning ``guesses'' at the answer.
We conclude with prospects and outlook in Section 7.

\section{Graph Laplacians}\setall

We commence with some rudiments from the theory of graphs, especially in relation to the Laplacian and its analogue for manifolds.
A graph $G$ is a pair $(V,E)$ where $V$ is a set of {\bf vertices} (we will also refer to them as nodes interchangeably) and $E$, a set of {\bf edges} (or arrows, whenever the edges are directed).
$V$ could be infinite but in this paper we will consider only finite graphs where the vertices in $V$ can be labeled as $v_{i = 1, \ldots, n}$.
A directed arrow $i \to j$ links $v_i$ to $v_j$ and a self-adjoining arrow $i \to i$ is called a loop.
A path in $G$ is a sequence of arrows $\{ v_{i_0} \to v_{i_1}, \ v_{i_1} \to v_{i_2}, \ v_{i_2} \to v_{i_3} , \ldots ,
v_{i_k} \to v_{i_{k+1}} \}$.
If the start and ending points of a path is the same node, i.e., $v_{i_0} = v_{i_{k+1}}$, then the path is called a {\bf cycle}.

As defined, we are allowing for multiple arrows between nodes. 
When we disallow (1) multiple edges  between nodes and (2) loops which link a node to itself, then $G$ is called {\bf simple}. 
Note that simple graphs {\em do allow} for cycles.
Furthermore, if we  ignore orientation of all arrows (and refer to them simply as edges), we have undirected graphs.
In congruence with the data available to us, which we will discuss shortly, we henceforth focus on finite, simple, undirected graphs, unless otherwise stated.

We now define some standard concepts in the theory of graphs.
\begin{definition}
The {\bf adjacency} matrix of $G = (V,E)$ with $|V| = n$ vertices is an $n \times n$ matrix $A_{ij} = \left\{
\begin{array}{lcr}
1 \ , & & \mbox{ if edge } i \to j \in E \\
0 \ , & & \mbox{otherwise}. 
\end{array}
\right.$

We use $\sim$ for adjacency: $v_i \sim v_j$ means nodes $i$ and $j$ is linked by an edge $i \to j$.
The {\bf degree} of vertex $v_i$ is the number of its neighbours:
$d(v_i) = \sum\limits_{v_i \sim v_j} 1$.
We record the degrees into a diagonal matrix, the {\bf degree matrix} $D_{ii}$ whose $i$-th diagonal entry is the degree of vertex $v_i$.
\end{definition}
In the case of our undirected graphs, $A$ is symmetric.
Moreover, for simple graphs the diagonal entries are all 0 (no loops) and all non-zero entries are 1 (no multi-edges).

Thus prepared, we can define the graph Laplacian simply as
\begin{definition}
For graph $G = (V,E)$, the Laplacian is defined as
\[
L := D - A \ ;
\qquad
L_{ij} = \left\{
\begin{array}{lcl}
d(v_i) && i = j \\
-1 && i \neq j \mbox{ and $v_i \sim v_j$} \\
0 && \mbox{otherwise} \ .
\end{array}
\right.
\]
\end{definition}
We remark that  for our simple undirected graphs with $n$ nodes, this is a symmetric $n \times n$ matrix, with diagonal entries being the degrees of the nodes.

\subsection{Rudiments on the Laplacian}
As defined above, it may seem obscure as to why $L$ is called the Laplacian.
To facilitate our grasp, we require several concepts:
\begin{definition}
The {\bf incidence} matrix $\nabla$ of $G = (V,E)$ is an $|E| \times |V|$ matrix with subscript $i=1,\ldots,|V|$ and $e = 1, \ldots, |E|$, such that for each edge $i \to j$, 
$\nabla_{i e} = -1$, $\nabla_{j e} = 1$, and 0 otherwise.
In other words, $(-1,+1)$ records the initial and final vertex of each edge, indexed by the rows. 
\end{definition}
For undirected graphs, we can choose some arbitrary but fixed orientation and whereby define $\nabla$.

Next, we can define a (real-valued) function $f : G = (V,E) \to \IR$ on a graph by assigning a real value to each node as $f(v_i)$.
We denote the space of such functions as
\begin{equation}
V^\IR := \{ f : V \to \IR \} \ .
\end{equation}
Subsequently, the adjacency and Laplacian matrices of $G$ can be seen as operators on functions on $G$.
For example, the adjacency matrix can be written as
\begin{equation}
A : f \longrightarrow g:=A(f) \ ; \ g(v_i) = \sum\limits_{i \sim j} f(v_j) \ .
\end{equation}
For our undirected graphs, this gives a convenient symmetric quadratic form $f^T A f = \sum\limits_{e \in E} f(v_i) f(v_j)$.
To illustrate, the following is a 4-vertex, 5-edged simple graph (the diamond graph), for which we write the adjacency matrix, choose an orientation, and write the incidence matrix:
\begin{equation}
\definecolor{wrwrwr}{rgb}{0.3803921568627451,0.3803921568627451,0.3803921568627451}
\definecolor{rvwvcq}{rgb}{0.08235294117647059,0.396078431372549,0.7529411764705882}
\begin{tikzpicture}[line cap=round,line join=round,>=triangle 45,x=0.8cm,y=0.8cm]
\draw [line width=2pt,color=wrwrwr] (0,2)-- (2,0);
\draw [line width=2pt,color=wrwrwr] (2,0)-- (0,-2);
\draw [line width=2pt,color=wrwrwr] (0,-2)-- (-2,0);
\draw [line width=2pt,color=wrwrwr] (-2,0)-- (0,2);
\draw [line width=2pt,color=wrwrwr] (-2,0)-- (2,0);
\begin{scriptsize}
\draw [fill=rvwvcq] (0,2) circle (2.5pt);
\draw[color=rvwvcq] (0.16,2.43) node {1};
\draw [fill=rvwvcq] (2,0) circle (2.5pt);
\draw[color=rvwvcq] (2.16,0.43) node {3};
\draw [fill=rvwvcq] (0,-2) circle (2.5pt);
\draw[color=rvwvcq] (0.16,-1.57) node {2};
\draw [fill=rvwvcq] (-2,0) circle (2.5pt);
\draw[color=rvwvcq] (-1.84,0.43) node {4};
\end{scriptsize}
\end{tikzpicture}
\begin{array}{l}
A = {\tiny \left(
\begin{array}{cccc}
 0 & 0 & 1 & 1 \\
 0 & 0 & 1 & 1 \\
 1 & 1 & 0 & 1 \\
 1 & 1 & 1 & 0 \\
\end{array}
\right)} \ ,  \\
D ={\tiny
\left(
\begin{array}{cccc}
 2 & 0 & 0 & 0 \\
 0 & 2 & 0 & 0 \\
 0 & 0 & 3 & 0 \\
 0 & 0 & 0 & 3 \\
\end{array}
\right)
}
\end{array}
\begin{tikzpicture}[line cap=round,line join=round,>=triangle 45,x=0.8cm,y=0.8cm]
\draw [->,line width=2pt,color=wrwrwr] (0,2) -- (2,0);
\draw [->,line width=2pt,color=wrwrwr] (0,2) -- (-2,0);
\draw [->,line width=2pt,color=wrwrwr] (-2,0) -- (0,-2);
\draw [->,line width=2pt,color=wrwrwr] (2,0) -- (0,-2);
\draw [->,line width=2pt,color=wrwrwr] (-2,0) -- (2,0);
\begin{scriptsize}
\draw [fill=rvwvcq] (0,2) circle (2.5pt);
\draw[color=rvwvcq] (0.16,2.43) node {1};
\draw [fill=rvwvcq] (2,0) circle (2.5pt);
\draw[color=rvwvcq] (2.48,0.09) node {3};
\draw [fill=rvwvcq] (0,-2) circle (2.5pt);
\draw[color=rvwvcq] (-0.04,-2.33) node {2};
\draw [fill=rvwvcq] (-2,0) circle (2.5pt);
\draw[color=rvwvcq] (-2.48,-0.01) node {4};
\end{scriptsize}
\end{tikzpicture}
\nabla = {\tiny
\left(
\begin{array}{cccc}
 -1 & 0 & 1 & 0 \\
 -1 & 0 & 0 & 1 \\
 0 & 0 & 1 & -1 \\
 0 & 1 & 0 & -1 \\
 0 & 1 & -1 & 0 \\
\end{array}
\right)
}
\end{equation}

As the notation suggests, $\nabla$ is a ``differential'' operator on functions on $G$: it is a co-boundary map in the sense that on an edge $i \to j$, $(\nabla f)(i \to j) = f(v_i) - f(v_j)$.
In analogy, we have
\begin{proposition}\label{prop:L}
The graph Laplacian is the ``square'' of the incidence matrix
\[
L = D - A = \nabla^T \nabla \ ;  \qquad  (L f)(v_i) = \sum\limits_{v_i \sim v_j} f(v_i) - f(v_j) \ ,
\]
irrespective of the choice of orientation.
\end{proposition}
\noindent {\em Sketch Proof:} To see this, we simply observe that $\sum_e \nabla_{i e} \nabla{e i}$ contributes 1 to each edge incident upon node $i$, which is then summed over the edges,  i.e., its gives the degree of node $i$. Similarly, $\sum_e \nabla_{i e} \nabla{e j}$ for $i \neq j$ contributes $-1 \cdot 1 = -1$ to the $(i,j)$-entry of $L$, which is the negative of the adjacency matrix.
Moreover, changing the direction of any edge in the direction assignment does not change the value of $\nabla_{i e} \nabla{e j}$ because one will be $+1$ and the other $-1$, whose product remains $-1$.

As we consider only undirected graphs, the relevant matrices $A$ and $L$ are symmetric, whereby giving us only real eigenvalues, these are respectively called the adjacency {\bf spectrum} and Laplacian spectrum of the graph.
Moreover, an important corollary of Proposition \ref{prop:L} is that the not only is the Laplacian spectrum real, its also non-negative:
\begin{equation}\label{specL}
\lambda_i = \mbox{Eigenvalues}(L) \ ; \quad 0 \leq \lambda_1 \leq \lambda_2 \leq \ldots \leq \lambda_{n = |V|} \ .
\end{equation}
This is seen by considering the eigensystem $L \phi = \lambda \phi$, which can be rewritten as $\lambda = \phi^T L \phi = \phi^T \nabla^T \nabla \phi = |\nabla \phi|^2 \geq 0$.

Finally, we remark that oftentimes we {\it weight} the edges of the graph by assigning a real postive number (the weight) $w_{ij}$ to each edge $i \to j$. Calling the set of weights $W$, we consider undirected weighted graphs $G = (V,E,W)$ and modify all of the above concepts accordingly.
The weighted adjacency matrix $A$ of $G = (V,E,W)$ is such that $A_{ij} = w_{ij}$ and 0 otherwise.
Likewise, the weighted Laplacian is
$L = D - A \ ;  (Lf)(v_i)= \sum\limits_{v_j \sim v_i} w_{ij}\big( f(v_i) - f(v_j) \big)$.

\subsection{Connection to Geometry}\label{s:geo}
The foregoing definitions and results are standard and can be found, for example, in \cite{FanBook}.
In the ensuing we will adhere to the discussions and conventions of \cite{ly}, in light of connections to differential and algebraic geometry.
First, we re-scale the Laplacian as
\begin{equation}\label{Lrw}
\Delta := D^{-1} L = I - D^{-1} A \ ; \qquad
(\Delta f)(v_i) = \frac{1}{d(v_i)}\sum\limits_{v_i \sim v_j} f(v_i) - f(v_j) \ .
\end{equation}
In the literature, this is often called the {\bf random walk normalized graph Laplacian} and is the one used in the programme of Yau et al. Likewise, we define a normalized version of the incidence (coboundary map):
\begin{equation}
\left| \nabla f (v_i) \right|^2 := \frac{1}{d(v_i)}\sum\limits_{v_i \sim v_j} \big( f(v_i) - f(v_j) \big)^2
\end{equation}
We make an important remark that in \cite{ly} (q.v., the first definition of $\Delta$ on p2), the Laplacian therein actually is the negative of  \eqref{Lrw}, which differs from some of the graph-theory literature.
However, this negative sign is compensated by its re-insertion in the eigen-equation on p4.
We will forego this double-negative and adhere to \eqref{Lrw}.

Now, upon defining a blinear operator $V^\IR \times V^\IR \to V^\IR$ as
\begin{equation}
\Gamma(f,g)(x) = \frac12 \big( \Delta( f g) - f \Delta g - g \Delta f \big) \ ,
\end{equation}
the concept of the curvature of a graph was introduced by \cite{fy}:
\begin{definition}\label{defRicci}
The {\bf Ricci curvature} of a graph is given by
\[
\Gamma_2(f,g)(x) = \frac12 \big( \Delta \Gamma( f g) - \Gamma(f, \Delta g) - \Gamma(g, \Delta f) \big)(x) \ .
\]
\end{definition}

We now recall from differential geometry that for a compact, smooth and complete Riemannian manifold $M$, there is {\bf Bochner's formula} which relates harmonic functions $u: \Delta u =0$ on $M$ to the Ricci curvature $Ric$ (cf.~\cite{cln}):
\begin{equation}
\frac12 \Delta ( \left| \nabla u \right|^2 ) = \left| \nabla^2  u \right|^2 + Ric(\nabla u, \nabla u) \ ,
\end{equation}
where $\nabla u$ is the gradient of $u$ and $Ric$, the Ricci curvature tensor, both with respective to the Riemannian metric on $M$.

One of the initial motivations of \cite{ly,ly2} was to have a discrete, graph-theoretic, versions of the above.
In particular, one has that \cite{ly2} (cf.~also \cite{glly})
\begin{theorem}[Lin-Yau]\label{thm:ly}
Let $G = (V,E)$ be a locally \footnote{
Note that the conditions here are more general than what we need and the graph itself can have an infinite number of vertices.
Locally finite means that at least  all vertex degrees are finite.
The supremum over all such degrees can, however, be infinite; hence we take the sup rather than max.
} finite graph and $d(G) := \sup\limits_{x \in V} d_x$ is the supremum over all vertex degrees (it can be that $d(G) = \infty$), then 
\[
\Gamma_2(f,f) \geq \frac12 (\Delta f)^2 + \left( \frac{1}{d(G)} - 1 \right) \Gamma(f,f) \ ,
\]
for any $f \in V^\IR$. 
\end{theorem}
We will return to address such inequalities in \S\ref{s:ineq}.

\section{The Wolfram Database of Simple Graphs}
\begin{figure}[!h!t!b]
\centerline{
\begin{tabular}{l}
(a1) \includegraphics[trim=0mm 0mm 0mm 0mm, clip, width=2in]{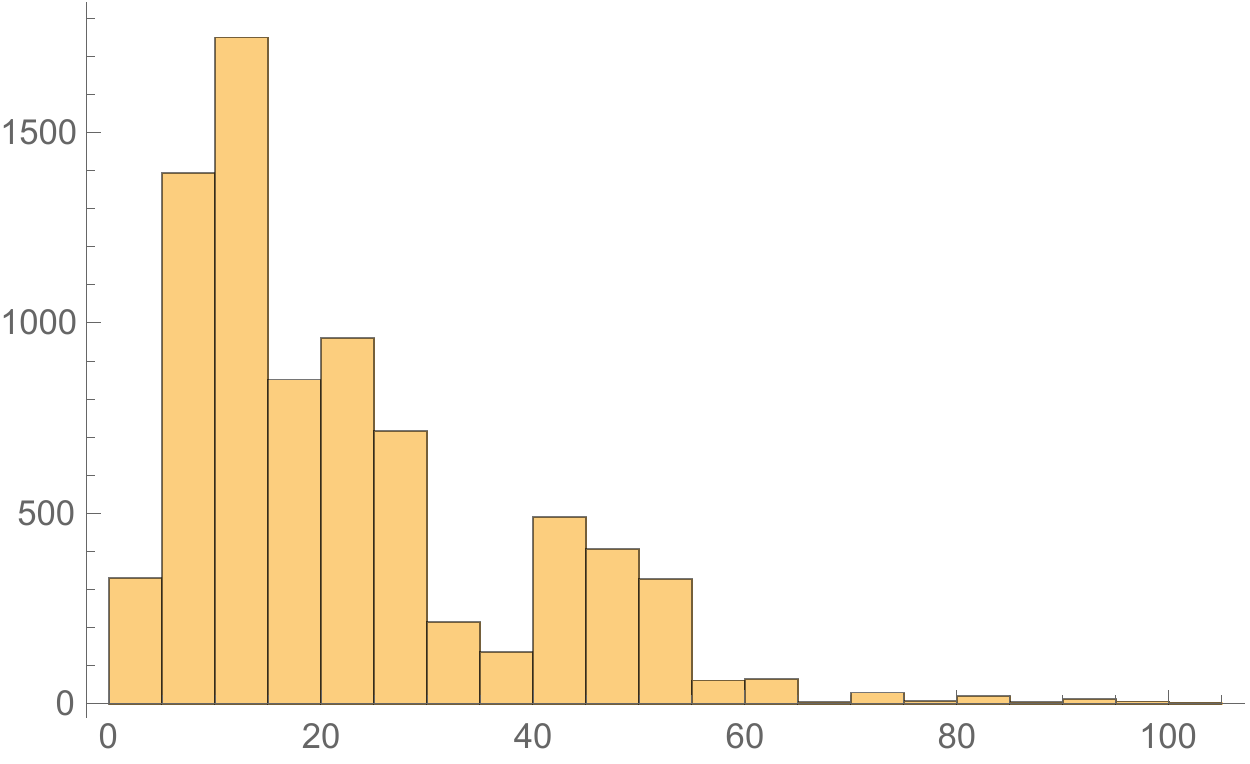}\\
(a2) \includegraphics[trim=0mm 0mm 0mm 0mm, clip, width=2in]{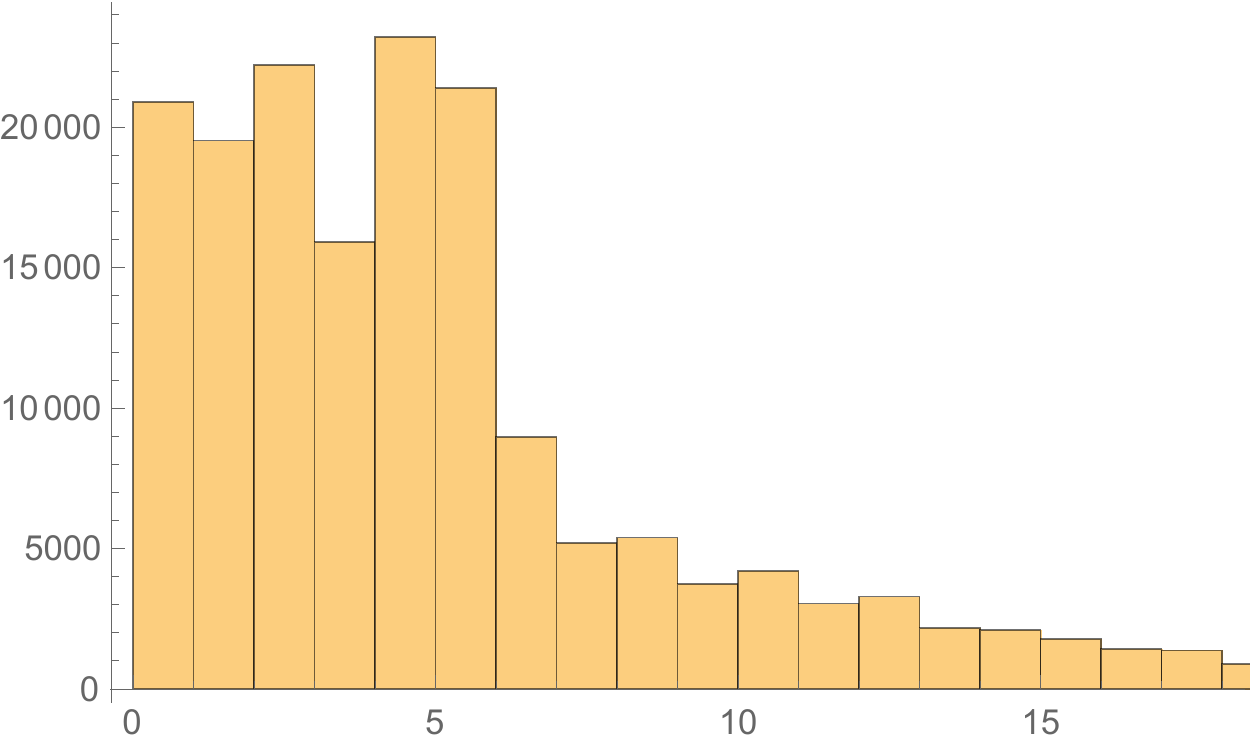}\\
\end{tabular}
(b)
\begin{tabular}{c}
\includegraphics[trim=0mm 0mm 0mm 0mm, clip, width=2.7in]{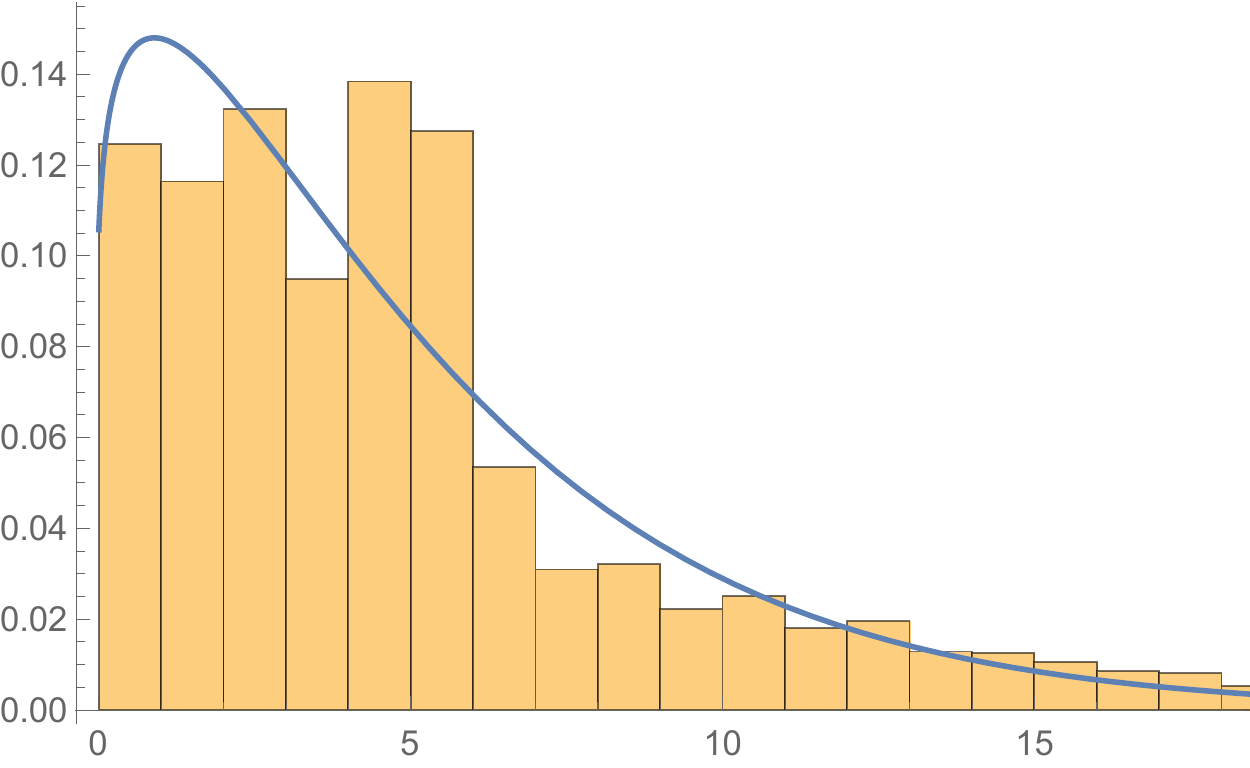}
\end{tabular}
}
\caption{{\sf {\small
(a1) The histogram of the number of vertices of the Wolfram database of 7785 undirected simple graphs up to 100 nodes;
(a2) The histogram of their Laplacian eigenvalues;
(b) The associated probability density (normalized histogram) and a fit by the  Weilbull distribution.
}}
\label{f:graphHist}}
\end{figure}

\setall
As a concrete play-ground, we take the graph database from Wolfram \cite{graphdata}, as implemented in Mathematica \cite{wolfram}, up to 100 vertices.
This is a list of undirected simple graphs, totaling 7785, which is a sizable set upon which we shall experiment.
We emphasize that this number is far less than the total number of known non-isomorphic, connected simple graphs up to 100 nodes, which proceeds exponentially \cite{numgraph} with the number of nodes as
\begin{equation}
\begin{array}{c}
1, 1, 1, 2, 6, 21, 112, 853, 11117, 261080, 11716571, 1006700565, 164059830476, \\
50335907869219, 29003487462848061, \ldots
\end{array}
\end{equation}
A histogram of the number of vertices (dimension of adjacency matrix) of the list of graphs is drawn in Part (a1) of Figure \ref{f:graphHist}

\begin{figure}[!!!t]
\centerline{
(a)
\includegraphics[trim=0mm 0mm 0mm 0mm, clip, width=1.5in]{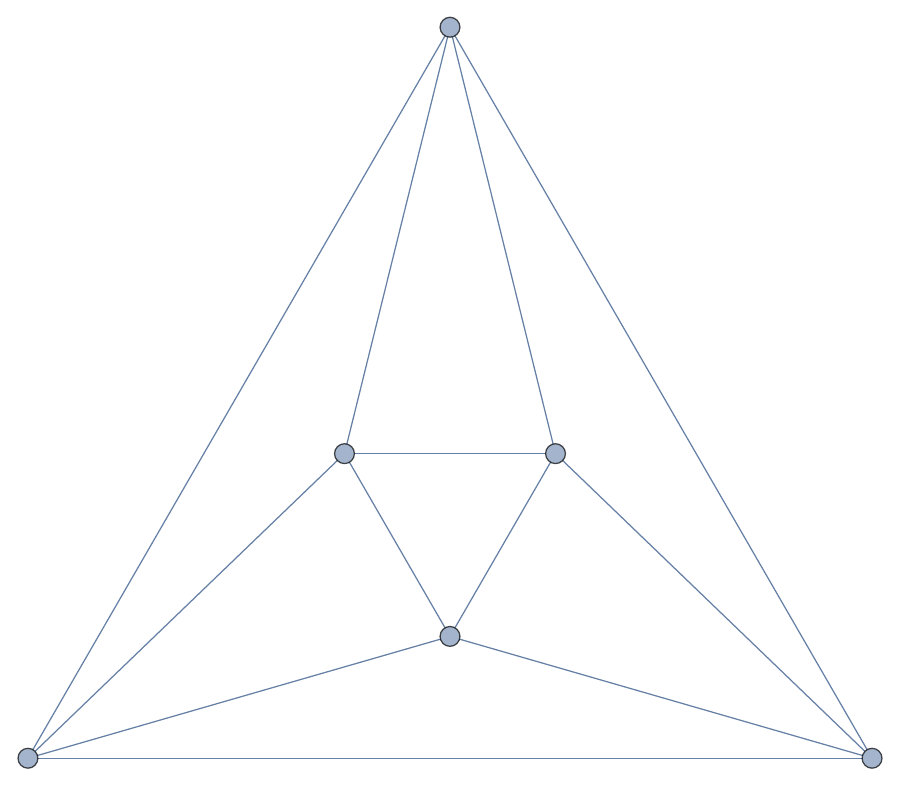}
(b)
\includegraphics[trim=0mm 0mm 0mm 0mm, clip, width=1in, angle=90]{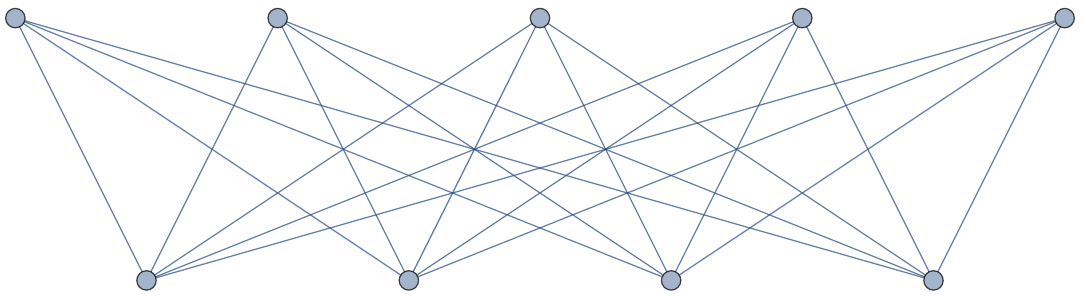}
}
\caption{{\sf {\small
(a) the Octahedral graph of 6 nodes; (b) the complete bipartite graph of 9 nodes from Wolfram {\sf GraphData[~]} database.
}}
\label{f:grapheg}}
\end{figure}

As explicit examples, we have Octahedral graph of 6 vertices in part (a) and the complete bipartite of 9 vertices in part (b) of
Figure \ref{f:grapheg}.
For reference, the Laplacians are $L_{\mbox{Octa}} =
{\tiny
\left(\arraycolsep=1.4pt\def\arraystretch{0.4}
\begin{array}{cccccc}
 4 & -1 & -1 & -1 & -1 & 0 \\
 -1 & 4 & -1 & -1 & 0 & -1 \\
 -1 & -1 & 4 & 0 & -1 & -1 \\
 -1 & -1 & 0 & 4 & -1 & -1 \\
 -1 & 0 & -1 & -1 & 4 & -1 \\
 0 & -1 & -1 & -1 & -1 & 4 \\
\end{array}
\right)}$ and $L_{\mbox{Bipartite} (5,4)} = {\tiny
\left(\arraycolsep=1.4pt\def\arraystretch{0.4}
\begin{array}{ccccccccc}
 5 & 0 & 0 & 0 & -1 & -1 & -1 & -1 & -1 \\
 0 & 5 & 0 & 0 & -1 & -1 & -1 & -1 & -1 \\
 0 & 0 & 5 & 0 & -1 & -1 & -1 & -1 & -1 \\
 0 & 0 & 0 & 5 & -1 & -1 & -1 & -1 & -1 \\
 -1 & -1 & -1 & -1 & 4 & 0 & 0 & 0 & 0 \\
 -1 & -1 & -1 & -1 & 0 & 4 & 0 & 0 & 0 \\
 -1 & -1 & -1 & -1 & 0 & 0 & 4 & 0 & 0 \\
 -1 & -1 & -1 & -1 & 0 & 0 & 0 & 4 & 0 \\
 -1 & -1 & -1 & -1 & 0 & 0 & 0 & 0 & 4 \\
\end{array}
\right) \ .
}$
We include in the Appendix a detailed walk-through of a particular example, exemplifying all the concepts encountered in the main body.

As mentioned in \eqref{specL}, all eigenvalues of the Laplacians of undirected graphs are non-negative.
To give an idea of the distribution of the spectrum, we show the histogram in part (a) of Figure \ref{f:graphHist}.
In part (b) of the figure, we find that a good fit is a Weibull distribution
\begin{equation}
p(x) = (x-\mu )^{\alpha -1} e^{-\left(\frac{x-\mu }{\beta }\right)^{\alpha }} \Theta(x-\mu) \ ; \qquad
\alpha = 1.16, \ \beta =  5.20, \  \mu = -0.05 \ ,
\end{equation}
where $\Theta(x)$ is the step function which is 1 for $x \geq 0$ and 0 for $x < 0$.
These should be compared with theoretical results on spectrum distributions of random graphs in \cite{gm,dj}
as well as that of the Laplacian eigenvalue distribution  \cite{dmt}.

\subsection{Preliminary Machine-Learning of Graph Properties}
With our database it is expedient that we  test whether some preliminary key graph-theoretical quantities can be machine-learned.
This is in the spirit of the citations in the introduction of whether various relevant quantities in string theory and mathematical physics can be deep-learnt.
In particular, the question of whether fundamental mathematical structures can be detected by AI in a supervised way has been pursued in the programme of \cite{He:2017aed,He:2019vsj,He:2019nzx,Alessandretti:2019jbs,He:2020eva}, ranging from number theory to representation theory, to computational geometry. Importantly, the idea is to use neural networks and classifiers that have {\it no prior knowledge} of the mathematics: indeed, with specifically chosen architectures guided by human intuition one achieves high precision (and even generate exact formulae as in \cite{Brodie:2019dfx}), but can generic networks find unseen patterns?
This seems to be the case in many of the aforementioned instances.

Our paradigm is clear.
Take our dataset $\cD$ of the Wolfram graphs.
Many important quantities have been computed and compiled. Take property $\cP$ and split $\cD = \cT \sqcup \cV$, the training set and validation set.
The training set then consists of association rules (labels) of the form $\cT = \{ A_i \to \cP_i \}$ with adjacency matrix $A_i$ for the $i$-th graph is to be associated with property $\cP_i$; an appropriate neural network or similar machine-learning algorithm is then trained with $\cT$.
The result is then validated against $\cV$ which the machine has {\em not} seen before in order to avoid over-fitting.
The validation is done, in the same spirit as a regression, by letting the algorithm {\it predict} the property $\cP$ for elements in $\cV$, after which we can compute some measure of ``goodness of fit''.
What we have described above is the paradigm of {\bf supervised learning}, where clear properties (labels) are assigned to the training set.

Typically, if $\cP$ is continuous, one can find a linear regression on the predicted versus actual values of $\cP$ on the validation set $\cV$, the fit should be as close to the line $y = x$ as possible.
If $\cP$ is discrete (and of a finite set of values, which can be labeled as $1,2,\ldots,k$), this is called ``$k$-categorical data''
(when $k=2$, this is binary classification). 
Here, a {\bf confusion matrix} $C_{ij}$ may be set up, where $i$ indexes the actual $k$ categories and $j$ indexes the machine-predicted $k$-categories so that 1 is added to each such pair. In the ideal case of perfect prediction, $C_{ij}$ is a diagonal matrix.

Several standard measure of the quality of the prediction (q.v.~e.g., \cite{sheskin}).
There is the naive {\bf precision} which is the percentage of correct predictions:
\begin{equation}
P_N := \sum_i C_{ii} \bigg/ \sum_{i,j} C_{ij} \ .
\end{equation}
The naive precision is too crude in that it does not distinguish amongst the diagonal entries.
A good measure (cf.~\cite{matthew}) for $k$-categorical classification is the {\bf Matthews correlation coefficient}:
\begin{equation}
\phi := {\frac {\sum\limits_{k}\sum\limits_{l}\sum\limits_{m}C_{kk}C_{lm}-C_{kl}C_{mk}}{{\sqrt {\sum\limits_{k}(\sum\limits_{l}C_{kl})(\sum\limits_{k'|k'\neq k}\sum\limits_{l'}C_{k'l'})}}{\sqrt {\sum\limits_{k}(\sum\limits_{l}C_{lk})(\sum\limits_{k'|k'\neq k}\sum\limits_{l'}C_{l'k'})}}}}
\end{equation}
The value of $\phi \in [-1,1]$ where 1 means perfect matching, 0 means random correlation and $-1$ means anti-correlation.

In the case of binary classifications where there are only 2 categories, the confusion matrix is usually denoted as 
\begin{equation}
C = \mbox{
\begin{tabular}{cc|c|c|}
	\cline{3-4} 
	&  & \multicolumn{2}{c|}{{Actual}} \\ \cline{3-4}
	& &  True  (1) & False (0)\\ \hline
	\multicolumn{1}{|c|}{{Predicted}} & True (1) & True Positive
	($tp$) & False Positive ($fp$) \\  \cline{2-4}
	\multicolumn{1}{|c|}{{Classification}} & False (0) & False Negative
	($fn$)& True Negative ($tn$)\\ \hline
\end{tabular}
}
\end{equation}
and the Matthews correlation coefficient reduces
\begin{equation}
\phi =\frac{tp \cdot tn - fp \cdot fn}{\sqrt{(tp+fp)(tp+fn)(tn+fp)(tn+fn)}} \ .
\end{equation}
In such binary classifications, another measure is the F1-score, which also needs to be close to 1 for a good prediction:
\begin{equation}
F_1 :=\frac{2}{ \frac{1}{TPR} + \frac{1}{Precision}} \ , \quad
	\begin{array}{cc}
	\text{TPR} := \frac{tp}{tp+fn} \ , \quad & \text{FPR} := \frac{fp}{fp+tn} \,,\\
	\text{Accuracy } p := \frac{tp+tn}{tp+tn+fp+fn} \ ,\quad &\text{Precision} := \frac{tp}{tp+fp}\,.
	\end{array}
\end{equation}
where TPR (FPR) stands for true (false) positive rate.

Armed with the rudiments, we can immediately turn to some important properties as warm-up exercises to our machine-learning perspective on graphs.
Now, the intersection of graph theory and machine-learning is taking shape recently (q.v., e.g., \cite{cohen,nikolaos} and references therein).
However, we point out that our approach is paradigmatically different here.
We are not employing results from graph theory to establish and study neural networks.
Instead, we are doing the opposite of using standard ML technique to see whether crucial graph properties can be learnt by supervision.
This is in line with the programme of the first author over the last two years in trying to see whether mathematical structures can be machine-learnt (cf.~summary talk at \cite{hetalk}).

In the following, we will perform what is called {\bf 5-fold cross-validation}, where we split (sequentially is good enough) $\cD$ into 5 groups, train on a subgroup of these and validate against the remaining.
In particular, we will show performance with 
\begin{align}
\nn
&\mbox{Training set } := 
\mbox{group} \{1,2,3,4,5\} \setminus \{i\} \ , \\
& \mbox{Validation set} :=
\mbox{group} \{i\} \ ; \qquad i = 1,\ldots,5 \ .
\end{align}
In other words, we train on 80\% of the data, and validate on the remaining 20\%, 5 different times. The resulting measures of fit are recorded so as to compute the average and error.
In general, one can perform $n$-fold cross-validation, but 5-fold is standard.
This way, we avoid both over-fitting, in that the validation set is never seen during training, as well as sample bias, in that we have trained and validated on various different combinations of the data.

Several technical comments are at hand. 
First, we are going up to graphs with 100 nodes, whose adjacency matrices have $10^4$ (albeit sparse) entries. This becomes computationally too intensive for a laptop computer.
Thus, for the following, we will take the first 5000 of the database (we still employ the full set for statistics and for some unsupervised ML later). 
As can be seen from Figure \ref{f:graphHist} (a1), the 
vast majority of the available graphs are at small number of vertices (this is purely for the easy of computation and compilation in the database). In fact, up to 5000, the largest is 25 nodes.
We subsequently pad all adjacency matrices of the graphs with zeros to the right and to the bottom so that all $A_i$ are $25 \times 25$.

Second, it is important that equivalent representation of the input be built-in.
Indeed, given a graph, any relabeling of the vertices is equivalent: this amounts to the {\it same} row/column permutation of the adjacency matrix \footnote{
This should be contrasted with Cayley multiplication table of finite groups, as done in \cite{He:2019nzx}, where {\it independent} row and column permutations are allowed.
}.
We typically perform 20 random permutations to each adjacency matrix $A$ and assign then the same property, where increasing the size of the data to $10^5$ (in practice, the number is smaller due to possible same random permutations).
As a technical aside, we perform the permutations {\em before} the zero-padding since this clearly makes more sense.

Finally, we will shuffle the data completely so that the sizes of the graphs are not ordered.
This is to avoid the bias of seeing only graphs of a certain number of vertices and validating against those of a different number.
Of course, one could purposefully do this so as to try to extrapolate to more complicated graphs from simpler ones, as was done in the spirit of \cite{Bull:2018uow}.
We leave discussions on this extrapolation to later.
For now, in the 5-fold cross validation, we will always be training on an assortment of graphs of varying sizes and complexity.

\subsubsection{A Graphical Miscellany}
Before we move on to the chief quantity of our concern, viz., the graph Laplacian.
Let us warm-up with machine-learning of some well-known properties of graphs
(again, q.v.~the Appendix for some explicit examples).
\paragraph{Planarity: }
One of the most important properties of a graph is whether it is {\bf planar}, in other words, whether the graph can be embedded into a plane and therefore can be drawn so that no edges cross except meeting at the nodes.
This a clear binary classification problem: adjacency matrix $A \to$ yes/no for planarity, well adapted to our supervised learning philosophy.

At our 5-fold validation, we find that
\begin{equation}
P_N \simeq 0.812 \pm 0.004, \ F_1  \simeq 0.832 \pm 0.004, \ \phi \simeq 0.619 \pm 0.009 \ ,
\end{equation}
using a logistic regression classifier (several other methods were tested but this simple regression seems to be the optimal).
In other words, having seen 20\% of the data, the ML has predicted (in under a minute on an ordinary laptop) whether a graph is planar by ``looking'' at the adjacency matrix, to  rather good confidence.

Now, there is a generalization of planarity, which is called {\bf graph skewness}.
This is the minimal number of edges to remove which would render the $G$ planar.
We remark that one could try to set up the complete prediction of the skewness.
However, the variation here is enormous within our dataset and ranges from 0 (planar) to many thousands.
We could instead make a 3-category classification: 0 (planar), 1 (skewness = 1) and 2 (skewness $>1$).
Doing so gives a slightly less significant result than the above, with $P_N$ and $\phi$ dropping down to 
$0.747 \pm 0.005$ and $0.597 \pm 0.008$ (note that $F_1$ is not defined for non-binary classifications).
Upon examining the confusion matrix, we see that the majority of mis-classifications is due to differentiating 1 and 2.
In other words, it is a little more difficult for the ML to distinguish precise skewness, and the quick decision on planarity is much easier.

Next, one could analyze the {\bf genus}, i.e., the genus of the Riemann surface onto which $G$ can be embedded. Clearly, genus $g=0$ means that the graph is planar.
Again, we can split the situation into  $g=0$ (planar), $g=1$ (doubly periodic) or $g > 1$.
This is in parallel to the trichotomy of Riemann surfaces, whether as complex varieties, they are Fano, Calabi-Yau, or general type, admitting, respectively, positive, zero, or zero curvature.
Here, we permute each adjacency matrix within the 3 categories by 10, 20 and 180 respectively in order to reach a more balanced data-set of around 16K each.
It turns out that this division is very much amenable to machine-learning, and we find at our 5-fold validation, that
\begin{equation}
P_N \simeq  0.814 \pm 0.003 \ , \quad \ \phi \simeq  0.721 \pm 0.005 \ ,
\end{equation}
which is even better than planarity recognition (again, there is no $F_1$-score to report since it is  a ternary classification here).

\paragraph{Chromatic Number: }
The minimal number of colours needed to colour the vertices such that no two adjacent vertices share the same colour is the chromatic number and is another important combinatorial quantity.
Here, like skewness, the chromatic number has a large variation within our dataset.
Thus, for starters, we can split this into a binary classification problem: whether the chromatic number $ch$ is, say, $>2$ or $\leq 2$.
We remark that this distinction of $ch \leq 2$ or not is important, because it is inequivalent to whether $G$ is {\bf bipartite} \footnote{i.e., whether all vertices of the graph can be split into two disjoint sets such that each edge connects two vertices, one from each such set; the study of bipartite graphs have become important in theoretical physics, from brane-tilings in AdS/CFT to SYM amplitudes.}.

Now, within $\cD$, those with $ch > 2$ dominates, so we will do more permutation enhancements for $ch \leq 2$ (we do 5 permutations for each $ch >2 $ and 30 for $ch \leq 2$, giving a more balanced set of around 30K for each category). 
At our 5-fold validation, then we find
\begin{equation}
P_N \simeq 0.773 \pm 0.005, \ F_1  \simeq 0.772 \pm 0.005, \ \phi \simeq 0.548 \pm 0.009 \ ,
\end{equation}
using a random forest classifier, which was found to be optimal.
This performance is, interestingly, worse than machine-learning planarity.

In light of the four-colour theorem, it is expedient 
\footnote{In terms of graphs, the statement is that the chromatic number of any planar graph is at most 4 (cf.~\cite{4colour}).} to focus on classification chromatic number for planar graphs.
This can be turned into a 3-category classification: chromatic number 2, 3, or 4 for a planar graph within our data set.
We enhance the data by permuting 30, 10 and 30 times respectively for the 2, 3, 4 chromatic categories, giving around 10K for each.
Trying on most standard classifiers, regressors, as well as neural networks with sigmoid activation functions, does not seem to produce results more significant than $\phi \simeq 0.5$.

\paragraph{Diameter: }
To give a notion of how big a graph is, one typically uses the {\bf diameter} of $G$: first, find the shortest distance between any pair of vertices (i.e., the minimal number of edges forming a path that is needed to go from one to the other), then, the diameter is the maximum amongst all these distances.
Since this has a large variation for our data-set, we can roughly divide the diameter $D(G)$ to be in 3 categories:
$D(G) \leq 2$, $3 \leq D(G) \leq 4$, and $D(G) > 4$.
Permuting the first two categories by 10 and the third by 30 to give about 30K per category, we have a roughly balanced 3-category classification problem.
Then, at our 5-fold validation, we find
\begin{equation}
P_N \simeq 0.765 \pm 0.004,  \ \phi \simeq 0.647 \pm 0.005 \ ,
\end{equation}
using a gradient boosted tree classifier, which was found to be optimal.
Again, the performance is quite good and the computation, under a minute per epoch of training.

\paragraph{Girth: }
Another important measure of the size of the graph is the {\bf girth},
\begin{equation}\label{girth}
\mbox{Girth}(G) := \mbox{minimum over the lengths of all cycles in } G; \quad \infty \mbox{ if $G$ is acyclic}.
\end{equation}
Within our dataset, the girth varies from 3 (around 4500 graphs), 4 (around 1400), and more than 4 (around 1400, including the acyclics).
Enhancing the data with random permutations for these 3 categories by 5, 30 and 40 respectively give a fairly balance set of around 15-20k for each of the 3 categories.
A classifier is then trained (optimized between nearest-neighbour and decision-tree) and our 5-fold validation gives
\begin{equation}
P_N \simeq 0.771 \pm 0.017,  \ \phi \simeq 0.656 \pm 0.026 \ ,
\end{equation}
in a matter of minutes.

As a parallel problem, one could consider the binary classification of
\begin{equation}
\mbox{Girth}(G) = \mbox{ or } \neq \infty \ .
\end{equation}
This is a fundamental characteristic of $G$ because it tells whether it is acyclic (possessing any cycles).
While there are nice algorithms such as topological sort which decides this in polynomial time, let us see how such a problem responds to ML.
We find our 5-fold validation to give an extremely good behaviour of
\begin{equation}
P_N \simeq 0.954 \pm 0.001, \ F_1  \simeq 0.955 \pm 0.001, \ \phi \simeq  0.912 \pm 0.002\ ,
\end{equation}
using a gradient-boost decision tree.

\paragraph{Special Cycles: }
Two classic problems concerning graphs are cycles which traverse all vertices and edges.
Indeed, if a cycle traverses all edges exactly once, it is an {\bf Eulerian} cycle \footnote{Indeed, Euler's 1736 translation of the K\"onigsberg bridge problem to the study of such cycles started the subject of graph theory.}.
We can see whether machine-learning can distinguish graphs which possess Eulerian cycles or not as a binary classification.
In our database of 5000 graphs, randomly permuting the negatives 5 times and the positives 30 times gives a roughly balanced set of around 17K cases each.
At our 5-fold validations, we find that a random forest classifier obtains
\begin{equation}
P_N \simeq 0.731 \pm 0.015, \ F_1  \simeq 0.721 \pm 0.026, \ \phi \simeq  0.473 \pm 0.024\ ,
\end{equation}
which is again rather good.

Similarly, a cycle which traverses all vertices exactly once is called a {\bf Hamiltonian} cycle.
Again, we can turn our dataset into a binary classification problem of whether machine-learning can tell which graphs have a Hamiltonian cycle by looking at the adjacency matrix.
It should be emphasized that this is known to be an NP-hard problem so stochastically learning a classifier is important.
Here, random permutation of the negatives 15 times and the positives 6 times gives a fairly balanced set of around 20K cases each.
Using a random forest classifier, we find that
\begin{equation}
P_N \simeq 0.781 \pm 0.008, \ F_1  \simeq 0.770 \pm 0.009, \ \phi \simeq  0.564 \pm 0.017\ ,
\end{equation}
which as encouraging as deciding the presence of Eulerian cycles.

\subsubsection{Maximal Laplacian Eigenvalue}\label{s:maxL}
Having warmed up with the exploration of a collage of graphical properties to see how well they respond to machine-learning - with many having encouraging results - we now turn to the object of our main concern, viz., the Laplacian.
We will investigate $L$ first before proceeding to the normalized Laplacian $\Delta$ in the next section.

\begin{figure}[!h!t!b]
\centerline{
(a) \includegraphics[trim=0mm 0mm 0mm 0mm, clip, width=3in]{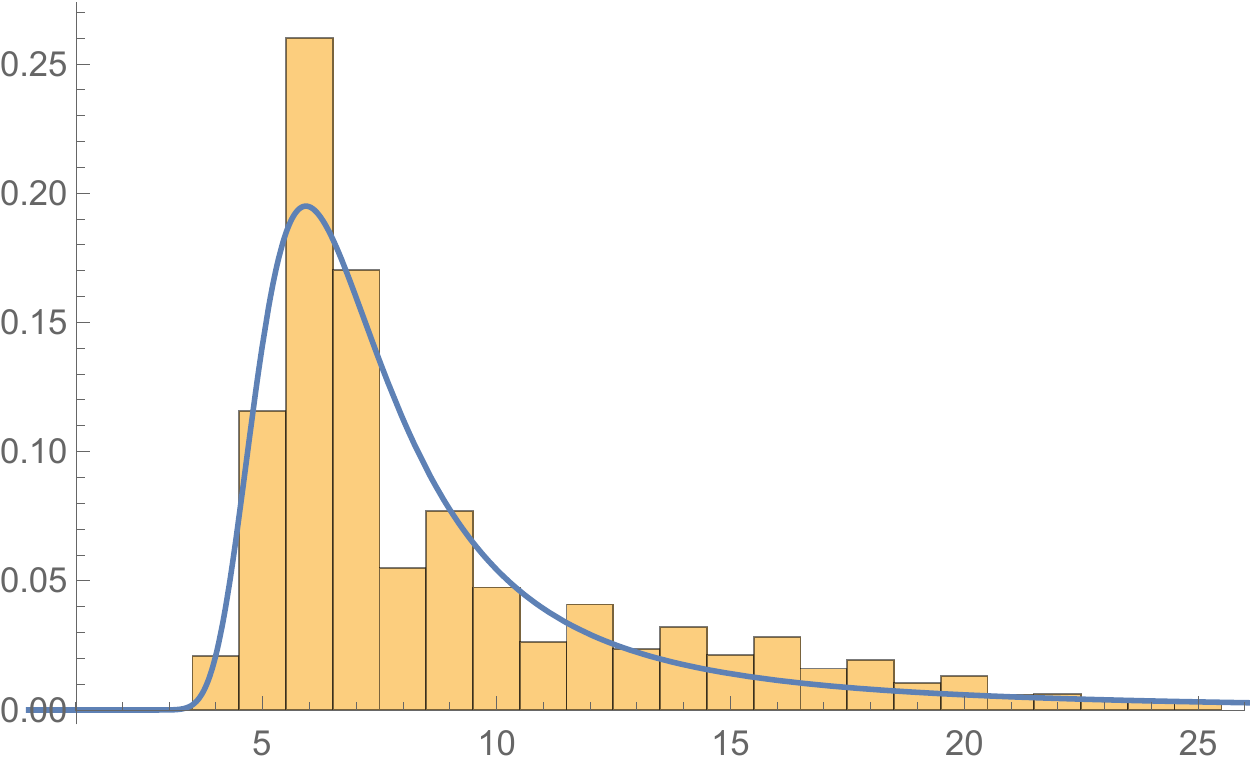}
(b)
\includegraphics[trim=0mm 0mm 0mm 0mm, clip, width=3in]{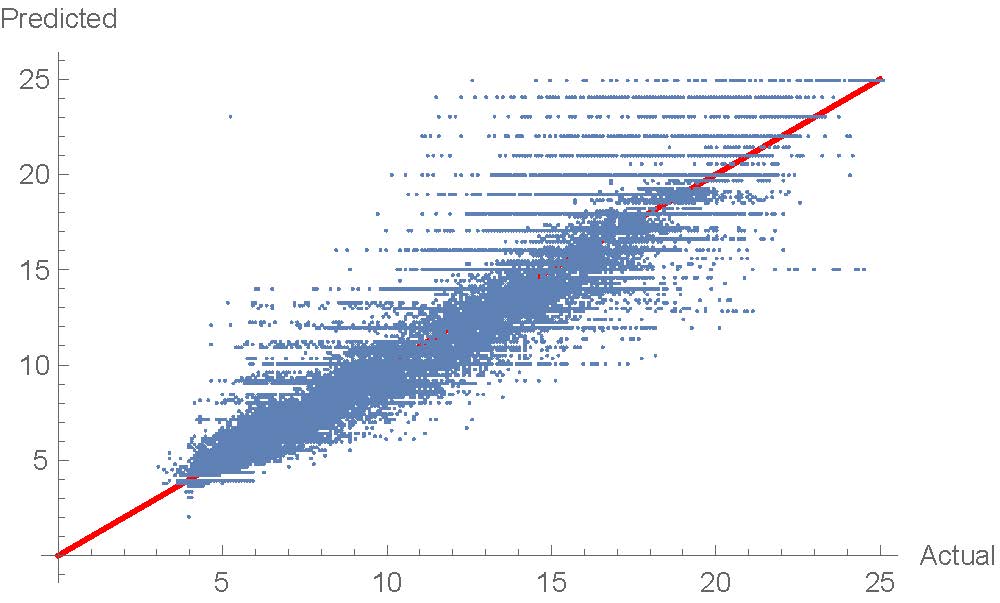}
}
\caption{{\sf {\small
(a) The histogram of the round-off of the maximal Laplacian eigenvalue for the first 5000 of our Wolfram database of undirected simple graphs.;
(b) The linear fit of the maximal Laplacian against that predicted by the NN in \eqref{NN} from the adjacency matrix.
}}
\label{f:maxEigenHist}}
\end{figure}

It is expedient to study the bounds on the Laplacian spectrum.
The lower, as shown earlier, is 0, so let us turn to the upper bound.
Within our sample of the first 5000 graphs, we first show a histogram, in Figure \ref{f:maxEigenHist}, of the {\bf maximal eigenvalue} of the Laplacian, rounded to the nearest integer.
The reason for this round-off will be explained momentarily.
A good fit of the distribution is found to be the Landau distribution:
\begin{equation}
p(x) = \int_0^\infty
\left(\frac{t}{\sigma }\right)^{-\frac{2 t}{\pi }} \sin (2 t) e^{\frac{t (\mu -x)}{\sigma }} dt \ , \qquad
\mu \simeq 6.204 \ , \quad
\sigma \simeq 1.456 \ .
\end{equation}

The reason we have taken the integer round-off is so that we can establish a discrete classification (we will move on to address the full, continuous problem shortly).
As can be seen from the distribution, there is a large variation here, consisting of all integers from 2 to 25.
In a matter of seconds on an ordinary laptop, using a logistic regression classifier for this 24-category problem, we find that our 5-fold validation gives
\begin{equation}
P_N \simeq 0.4982 \pm 0.007, \ \phi \simeq 0.415 \pm 0.007 \ .
\end{equation}
In other words, the regressor has predicted the correct values of the integer round of the maximal Laplacian eigenvalue to about 50\% with almost as much confidence.

Moving onto the full problem, we need to establish a predictor for maximal Laplacian spectrum in the form of
\begin{equation}
A_{ij} \longrightarrow \max\limits_{\lambda \mbox{ {\tiny Eigenvalue}}} \lambda(L) \ . 
\end{equation}
Of course, due to the high dimensionality of the input variables: there are $^{25}C_2 = 300$ degrees of freedom for the adjacency matrix $A_{ij}$, to find a simple function $f$ which, by some non-linear regression, that produces the maximal $\lambda$, would be rather difficult.

We establish a neural network (NN) in the form of a forward-propagating, 4 layer perceptron (MLP) which was found to very efficient in computing cohomology \cite{He:2017aed} and which has the form
\begin{equation}\label{NN}
\begin{array}{l}
\includegraphics[trim=0mm 0mm 0mm 0mm, clip, width=5in]{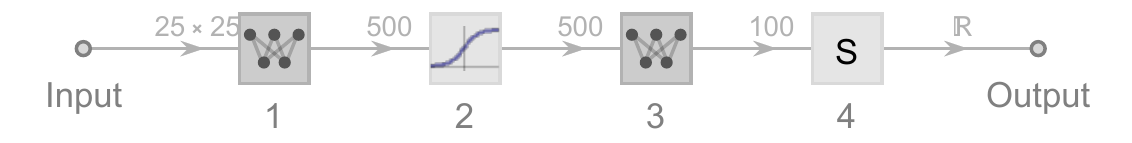}
\end{array}
\end{equation}
That is, the input layer is a $25 \times 25$ matrix, followed by a fully-connected linear layer taking it to 500 nodes, each of which is passed to a Sigmoid activation function $x \to (1 - e^{-x})^{-1}$. These 500 neural nodes are then mapped to 100 nodes in a fully connected layer, before finally summed to a real number, which should correspond to the output of the maximal $\lambda$.

Using an ADAM optimizer with batch size 64, and with supervised training of the above NN on the 20\% of the data, we can plot the predicted versus the actual maximal $\lambda$ for the unseen validation 80\%.
This is shown in Part (b) of Figure \ref{f:maxEigenHist}.
The scatter plot should be as close to the line $y=x$ (drawn in red) as possible and we do see indeed the predicted maximal $\lambda$ are so.
A best fit shows that
\begin{equation}\label{fit20maxEigen}
y  = -0.0619868 + 1.01097 x \ , \qquad R^2 = 0.93 \ ;
\end{equation}
the closeness of $R^2$, the coefficient of determination, to 1, shows that this is a good fit.

\begin{figure}[!ht!b]
\centerline{
\includegraphics[trim=12mm 0mm 0mm 0mm, clip, width=5in]{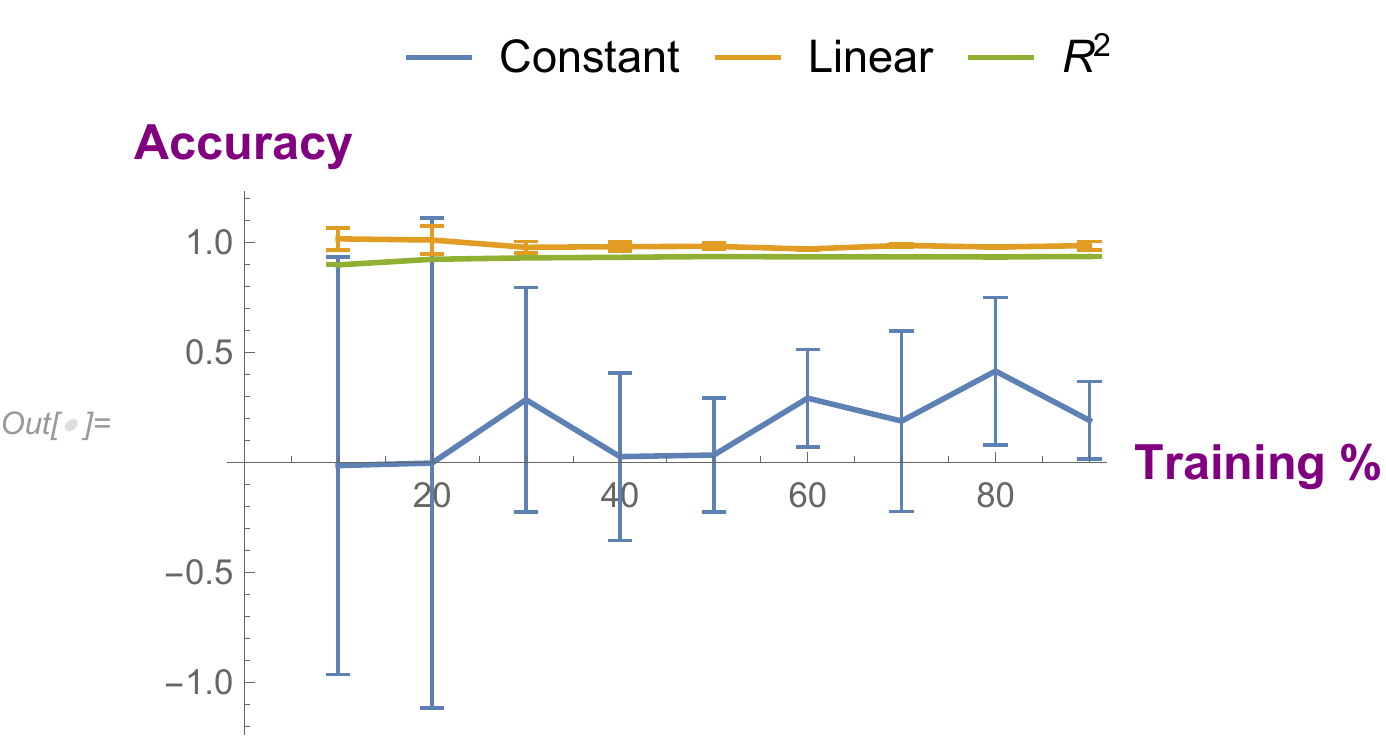}
}
\caption{{\sf {\small
The training curve for the neural network in \eqref{NN} for ML of the maximum Laplacian eigenvalue of our graph dataset.
A random sample of $x\%$ is trained and then the NN is validated on the remaining $(100-x)\%$; we let x range from 10 to 90 in increments of 10.
We show 3 measures of goodness of fit: the constant and linear terms of the linear regression of predicted and actual values, as well as the coefficient of determination $R^2$.
}}
\label{f:maxEigenTC}}
\end{figure}

We can get a glimpse of the machine-learning behaviour also by training on randomly chosen $x \%$ and validating on the complementary $(100 - x) \%$.
Plotting the goodness of fit of the ML gives us a {\bf training curve}.
This is is shown in Figure \ref{f:maxEigenTC}.
Training is done from 10\% to 90\% random sample from the graph database, and validated on the complement.
As in the aforegoing discussions, there is a list of predicted maximal Laplacian eigenvalues and a list of actual values, we perform linear regression $y = a x + c$ on these to obtain (i) the constant term $c$, (ii) the linear coefficient $a$, as well as (iii) the coefficient of determination $R^2$.
The value of $c$ should be close to 0, $a$ should be close to 1 and $R^2$, to 1.
We see that starting from a mere 20\%, the behaviour is already very good.
The error bars are collected from 5 rounds (epochs) of training: both $a$ and $R^2$ behave well from the beginning, whilst $c$ fluctuates less and less as we increasing training size.

\subsubsection{Spectral Gap}\label{s:gap}
On the other extreme, we can study the lowest Laplacian eigenvalue.
Of course, by Eq.~\eqref{specL}, the lowest eigenvalue is 0.
In fact, we have that (cf.~\cite{FanBook})
\begin{proposition}\label{0eigen}
The dimension of the nullspace of the graph Laplacian is equal to the number of connected components of $G$.
\end{proposition}
Thus, for all our graphs, there is {\em exactly one} 0-eigenvalue.
The next smallest eigenvalue (possibly with multiplicity) is called the (Laplacian) {\bf spectral gap}; it is therefore the first positive eigenvalue of $L$.

\begin{figure}[!h!t!b]
\centerline{
(a) \includegraphics[trim=0mm 0mm 0mm 0mm, clip, width=3in]{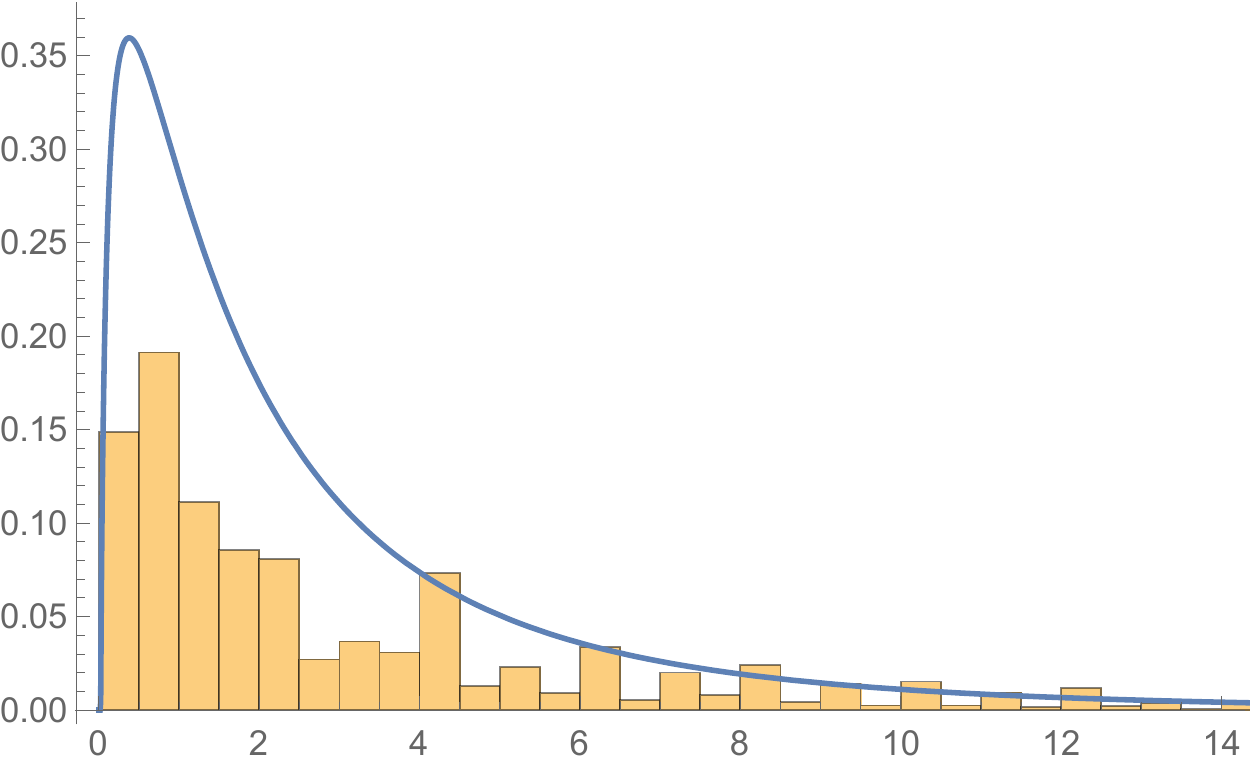}
(b)
\includegraphics[trim=0mm 0mm 0mm 0mm, clip, width=3in]{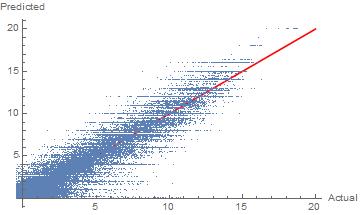}
}
\caption{{\sf {\small
(a) The histogram of the spectral gap (smallest positive Laplacian eigenvalue) for the first 5000 of our Wolfram database of undirected simple graphs;
(b) its linear fit against that predicted by the NN in \eqref{NN} from the adjacency matrix.
}}
\label{f:specGapHist}}
\end{figure}

We show the histogram  (normalized to a probability) in part (a) of Figure \ref{f:specGapHist}, together with its best fit as a probability distribution, which is found to be a generalized Gamma distribution
\begin{align}
\nn
p(x) &= (x-\mu)^{\alpha \gamma-1} \exp(-((x-\mu)/\beta])^\gamma) \ ; \qquad x > \mu \ , 
\\
&
\alpha \simeq 5.164, \
\beta \simeq 0.031, \
\gamma \simeq 0.387, \
\mu \simeq 0.023 \ .
\end{align}

As with the maximal eigenvalue, we perform ML by the NN in \eqref{NN}.
At 20\% random sample for training and validating against the complement 80\%, we can find the best fit line to be
$y \simeq 0.910 x + 0.700$ with $R^2 \simeq 0.897$, which is shown in Part (b) of Figure \ref{f:specGapHist}.
Interestingly, the fit here, whilst still rather good, is less well-bahaved than \eqref{fit20maxEigen}.

\section{Laplacian Spectra and Inequalities}\label{s:ineq}\setall
We have now practiced with ML on a plethora of graph theoretic quantities and can continue to study our protagonist the normalized Laplacian $\Delta$.
One of the classic results in differential geometry is the Li-Yau theorem on the bound of the eigenvalue of the Laplacian.
In particular, we have \cite{LiYau}
\begin{theorem}[Li-Yau]
Let $M$ be a compact Riemannian manifold of dimension $n$ and diameter
\footnote{
The diameter of a manifold is the supremum over all geodesic lengths on $M$.
}
 $\delta(M)$, with Ricci curvature bounded below by $(n-1)K$, then the first non-zero eigenvalue of the Laplacian is bounded below by
 \[
 \lambda \geq \frac{\exp \left[- (1 +  \sqrt{ 1 - 4(n-1)^2 \delta(M)^2 K}  \right]}{2(n-1)\delta(M)^2} \ .
 \]
\end{theorem}

As discussed in \S\ref{s:geo}, one of the motivations of studying graph Laplacians is to see whether there are such ``discrete'' analogues of these spectral bounds.
Now, it is a standard result in graph theory that the non-zero Laplacian eigenvalues are bounded by
(cf.~\cite{FanBook} and also q.v.~survey in \cite{mohar})
\begin{equation}\label{boundG}
\lambda \geq [\delta(G) V(G)]^{-1} \ ,
\end{equation}
where $\delta(G)$ is the diameter of the graph $G$ and $V(G) := \sum\limits_{v_i} d(v_i)$ is the volume of $G$, the sum over all degrees of all vertices.

This was improved by Lin-Yau \cite{ly,ly2} to be
\begin{equation}\label{boundLY}
\lambda \geq \left[ d(G) \delta(G) \exp\big( d(G) \delta(G) + 1) - 1 \right]^{-1} \ ,
\end{equation}
where, as was in Theorem \ref{thm:ly}, $d(G) := \sup\limits_{x \in V} d_x$ is the supremum over all vertex degrees.

\subsection{Eigenvalue Distributions of the Normalized Laplacian}
\begin{figure}[!h!t!b]
\centerline{
(a) \includegraphics[trim=0mm 0mm 0mm 0mm, clip, width=3in]{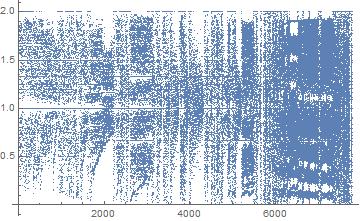}
(b) \includegraphics[trim=0mm 0mm 0mm 0mm, clip, width=3in]{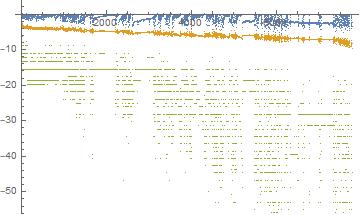}
}
\caption{{\sf {\small
(a) all eigenvalues for the 7785  connected simple graphs from the  Wolfram database, sequentially from left to right in the order of the database (roughly according to complexity);
(b) the minimal non-zero eigenvalue (spectral gap) of the the 7785 graphs, compared to the two bounds in \eqref{boundG} and \eqref{boundLY}.
}}
\label{f:eigenL}}
\end{figure}

To get an idea of the Laplacian spectrum, especially the lowest and the highest, for our database, we first order all the graphs from left to right in complexity (roughly the number of vertices and edges, as ordered in the Wolfram database) across the abscissa, then, on each point, we list all the Laplacian eigenvalues starting from the lowest to the highest (we recall from Proposition \ref{0eigen} that the lowest is a single 0 as all our graphs are connected).
This is shown in Part (a) of Figure \ref{f:eigenL}.
We see that all eigenvalues here are between 0 and 2.
The upper bound of 2 is a known result for the normalized Laplacian (equaling to 2 for bipartite graphs).

In Part (b) of the figure, we show a comparison between the magnitudes of (1) the smallest positive eigenvalue, the spectral gap $\lambda_{min}$, (2) the RHS of the standard bound \eqref{boundG}, and (3) the RHS of the Lin-Yau bound \eqref{boundLY}.
Due to the presence of the exponential, we show the log of the 3 quantities.
We see how $\lambda_{min}$ dominates, as is required.
Now, it might appear that the standard bound \eqref{boundG} is stronger than that of \eqref{boundLY}.
However, we need to remember that the latter applies to {\it locally} finite graphs which could have an infinite number of vertices and the subsequent sum could grow dramatically.

\subsection{Machine-Learning Spectral Bounds}
As in \S\ref{s:gap} for the non-normalized Laplacian $L$, we can now apply ML to the normalized Laplacian $\Delta$ of our concern, especially in light of the spectral inequalities of the Li-Yau/Lin-Yau type.

\subsubsection{Unsupervized Treatments}
We begin by applying some unsupervized techniques in order to better visualize the eigenvalue distribution.
The diagram in part (a) of Figure \ref{f:eigenL} is unenlightening. For one thing, the x-axis is an arbitrary ordering (part (b), of course, is meaningful in showing the comparison between the inequalities).
A more natural coordinate to take, for instance, is the ordered Laplacian eigenvalues.
To normalize, we left pad all matrices to the maximal dimension without dataset, viz, 100.
This way, we have a {\bf point cloud} $\Lambda \subset \IR^{100}_{\geq 0}$, each point of which is
\begin{equation}\label{eigen100}
\vec{\lambda} := (0, 0, \ldots, 0, \lambda_1, \lambda_2, \ldots, \lambda_{k})_{100} \ , \quad \lambda_{i-1} \leq \lambda_{i} \ ; \qquad \IR^{100} \supset \Lambda = \{\vec{\lambda}\} \ .
\end{equation}

\begin{figure}[!h!t!b]
\centerline{
(a) \includegraphics[trim=10mm 0mm 0mm 0mm, clip, width=3in]{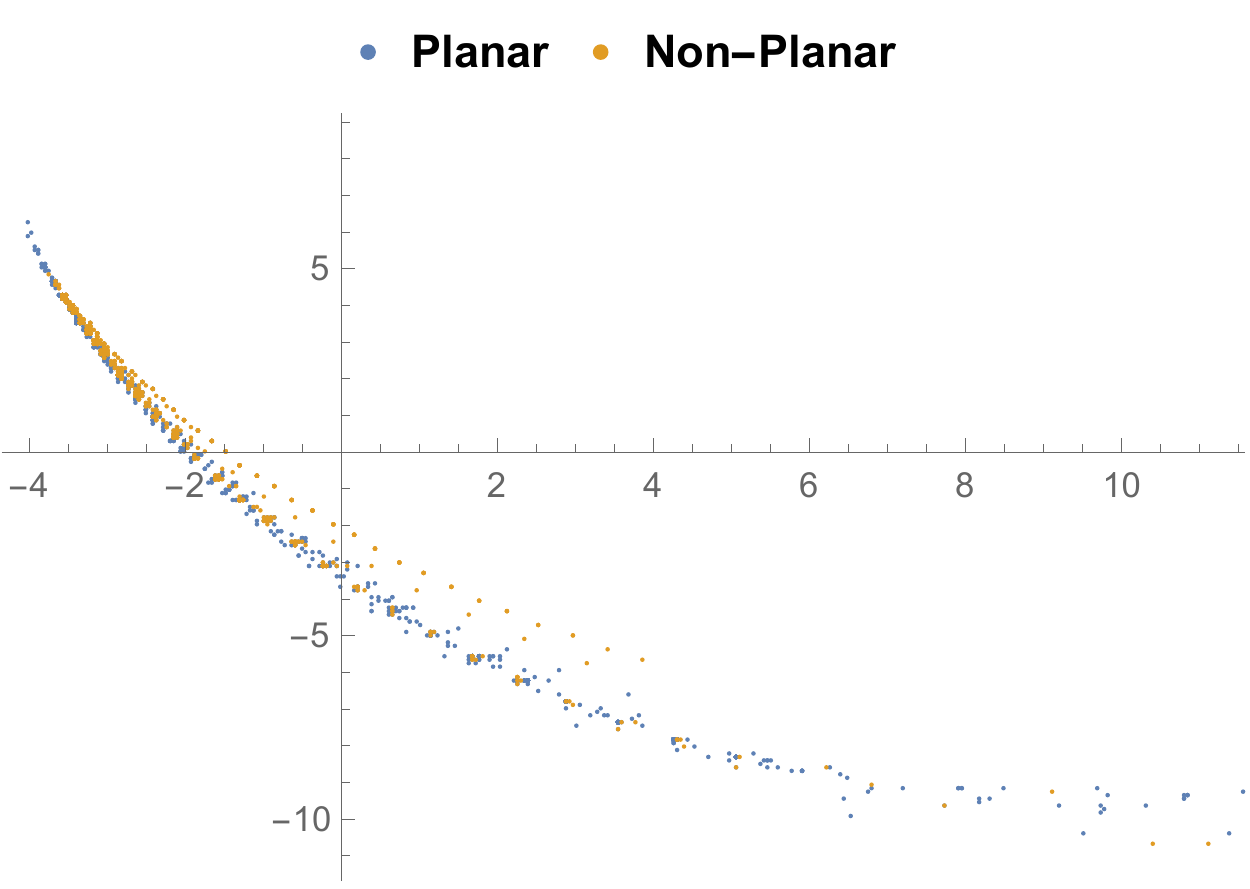}
(b) \includegraphics[trim=10mm 0mm 0mm 0mm, clip, width=3in]{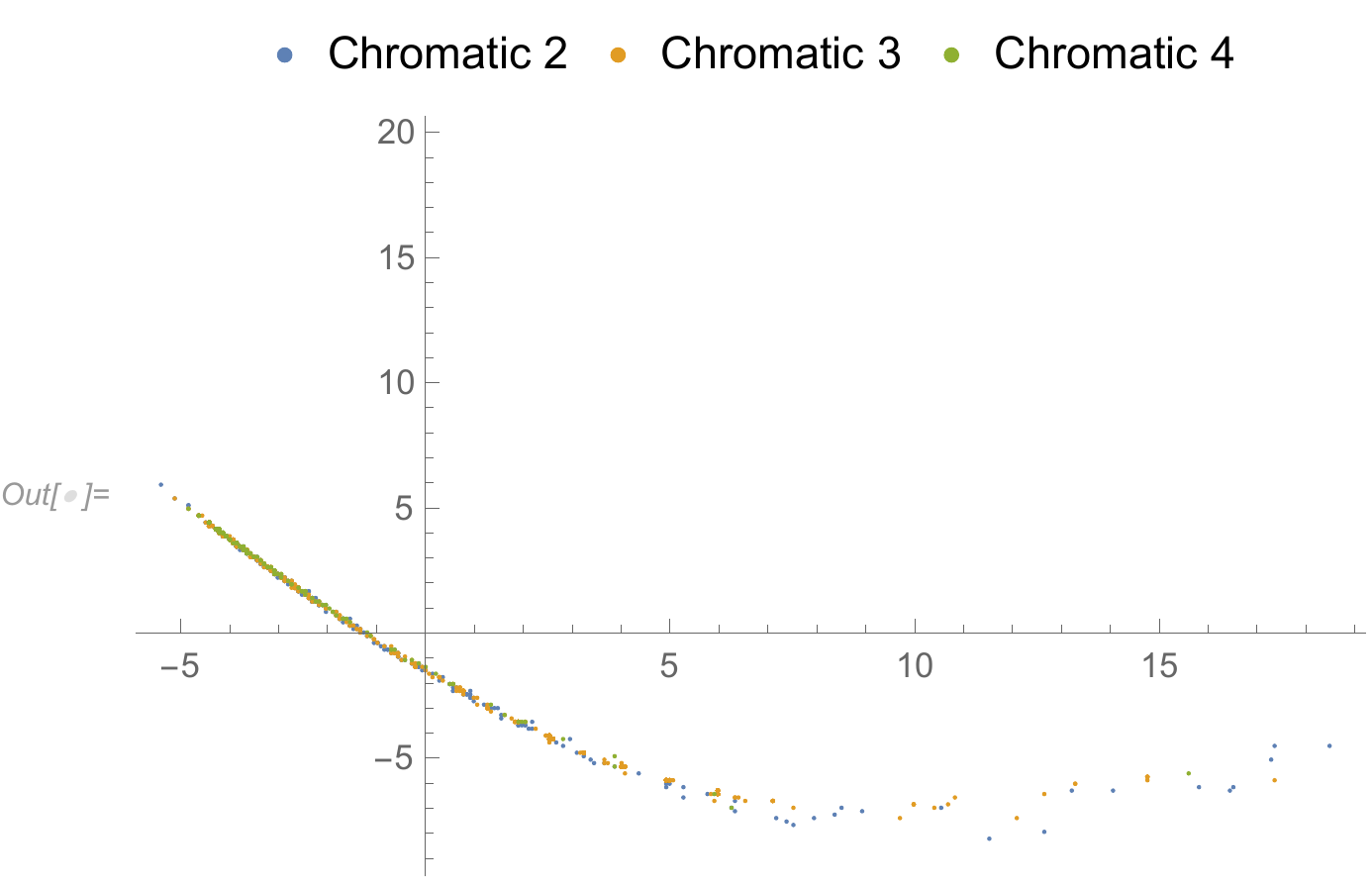}
}
\caption{{\sf {\small
(a) The PCA for our 7785 graphs, reducing the 100-dimensional point-cloud defined by the Laplacian spectrum as in \eqref{eigen100} to 2-dimension.
(b) Similarly, for the 2933 planar graphs, we plot the point-cloud, but distinguished by the 3-categories of chromatic numbers 2, 3, and 4.
}}
\label{f:PCAeigen}}
\end{figure}

This point cloud is of enormous dimension, but luckily unsupervized ML provides precisely the technique to visualize it: via the method of {\bf principle component analysis} (PCA). We can map $\IR^{100}$ to, for example, $\IR^2$.
The dimensionality reduction is shown in part (a) Figure \ref{f:PCAeigen}; it is an interesting semi-crescent shape.
As an illustration, we have separated the planar versus the non-planar graphs.
It is curious that the Laplacian spectrum of the non-planar graphs lie a little above that of the planar ones.
In part (b) of the same figure, we isolated the planar graphs (of which there are 2933). 
For these, the chromatic number can be 2, 3, or 4 due to the 4-colour theorem.
We perform a similar $(100 \to 2)$-dimensional reduction, and mark the 3 different chromatic numbers.
Here we still obtain a crescent shape but the 3 categories of chromatic numbers do not seem that much different.

\begin{figure}[!h!t!b]
\centerline{
(a) \includegraphics[trim=0mm 0mm 45mm 0mm, clip, width=2.5in]{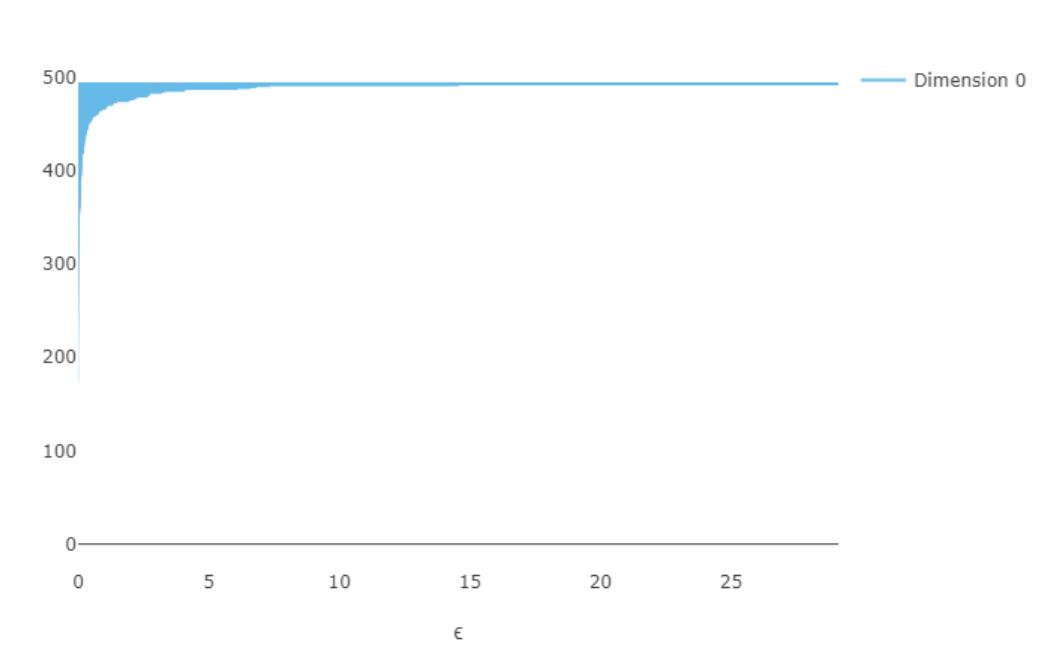}
(b) \includegraphics[trim=0mm 0mm 45mm 0mm, clip, width=2.5in]{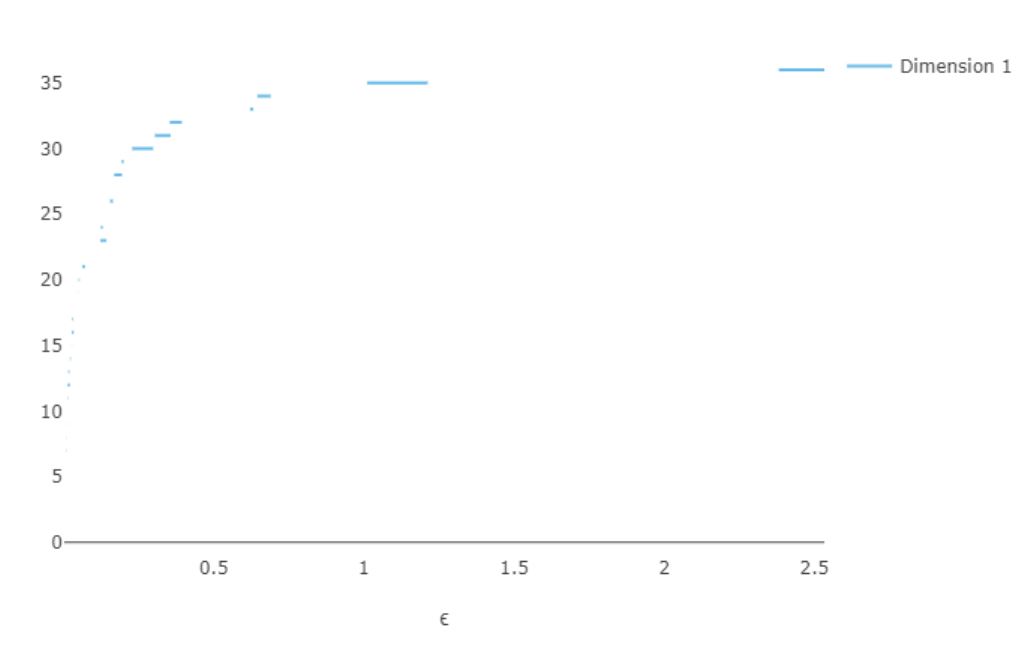}
}
\caption{{\sf {\small
The barcodes for 500 random samples of the 100-dimensional point cloud of the Laplacian spectrum of our graph dataset, indicating the persistent homology in dimension 0 and 1 (respectively parts (a) and (b)).
}}
\label{f:TDA}}
\end{figure}

Another way to visualize high dimensional data is to use topological data analysis (TDA) \cite{czcg} (cf.~recent review in \cite{tda}).
Here, we compute the persistent homology of the point cloud and present them in terms of so-called {\bf barcodes}.
Using the {\em Julia} implementation of {\em Eirene} from \cite{eirene}, we compute the zeroth and first persistent homology for the point cloud $\Lambda$ in \eqref{eigen100}; these are shown in parts (a) and (b) respectively in Figure \ref{f:TDA}. Note that computing barcodes is extremely expensive so around 7000 points in $\IR^{100}$ would be quite prohibitive; we subsequently took 500 random sample points to give an idea of the persistence.

\subsubsection{Supervized Learning}
Let us now perform supervised ML on the Laplacian eigenvalue.
First, we repeat the analysis in \S\ref{s:gap}, but now for the normalized Laplacian $\Delta$ (rather than the ordinary $L$ which was done there).
In other words, we have a labeled problem of
\begin{equation}
A_{ij} \longrightarrow \min\limits_{0 \leq \lambda \mbox{ {\tiny Eigenvalue}}} \lambda(\Delta) \ , 
\end{equation}
from the adjacency matrix to the spectral gap.
Suppose the ML has seen only 20 \% of the data at random, how does it predict the gap for the remaining 80\% ?

\begin{figure}[!ht!b]
\centerline{
(a)
\includegraphics[trim=10mm 0mm 0mm 0mm, clip, width=2in]{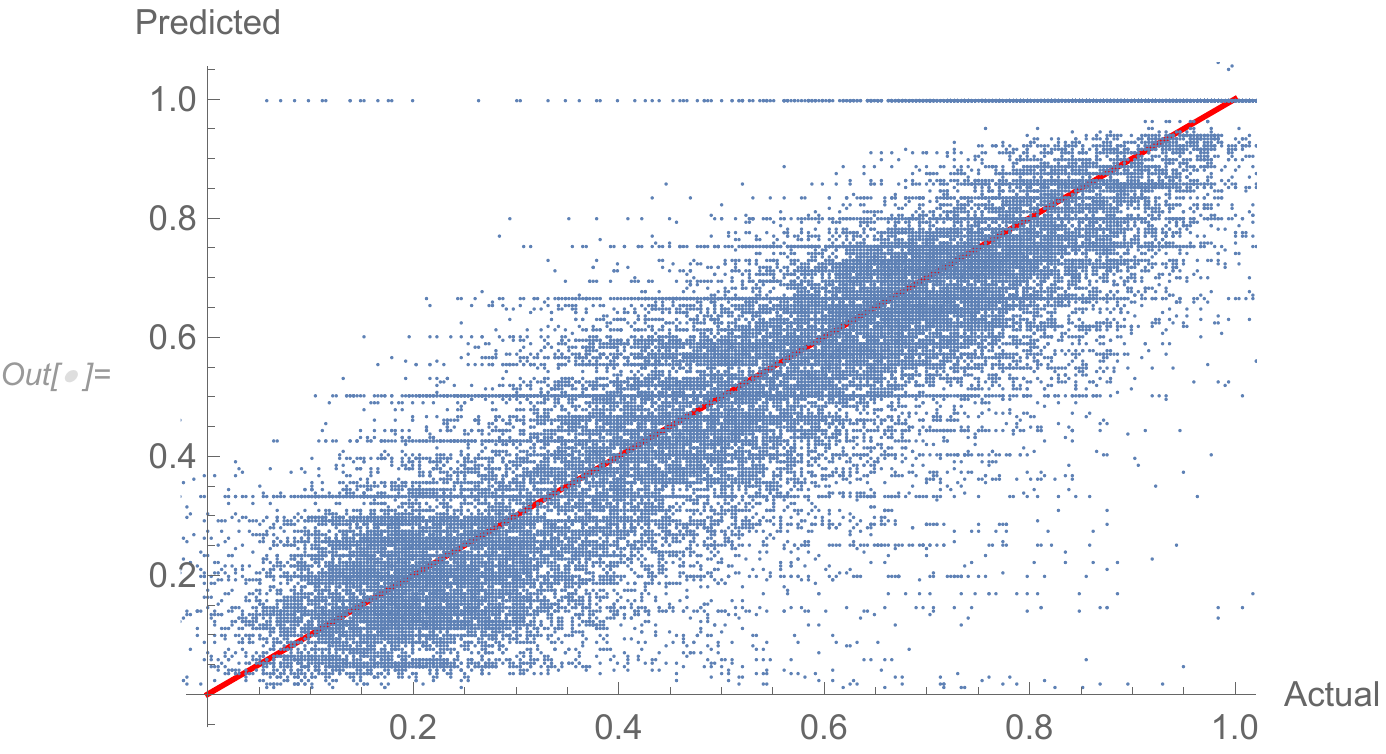}
(b)
\includegraphics[trim=10mm 0mm 0mm 0mm, clip, width=3.5in]{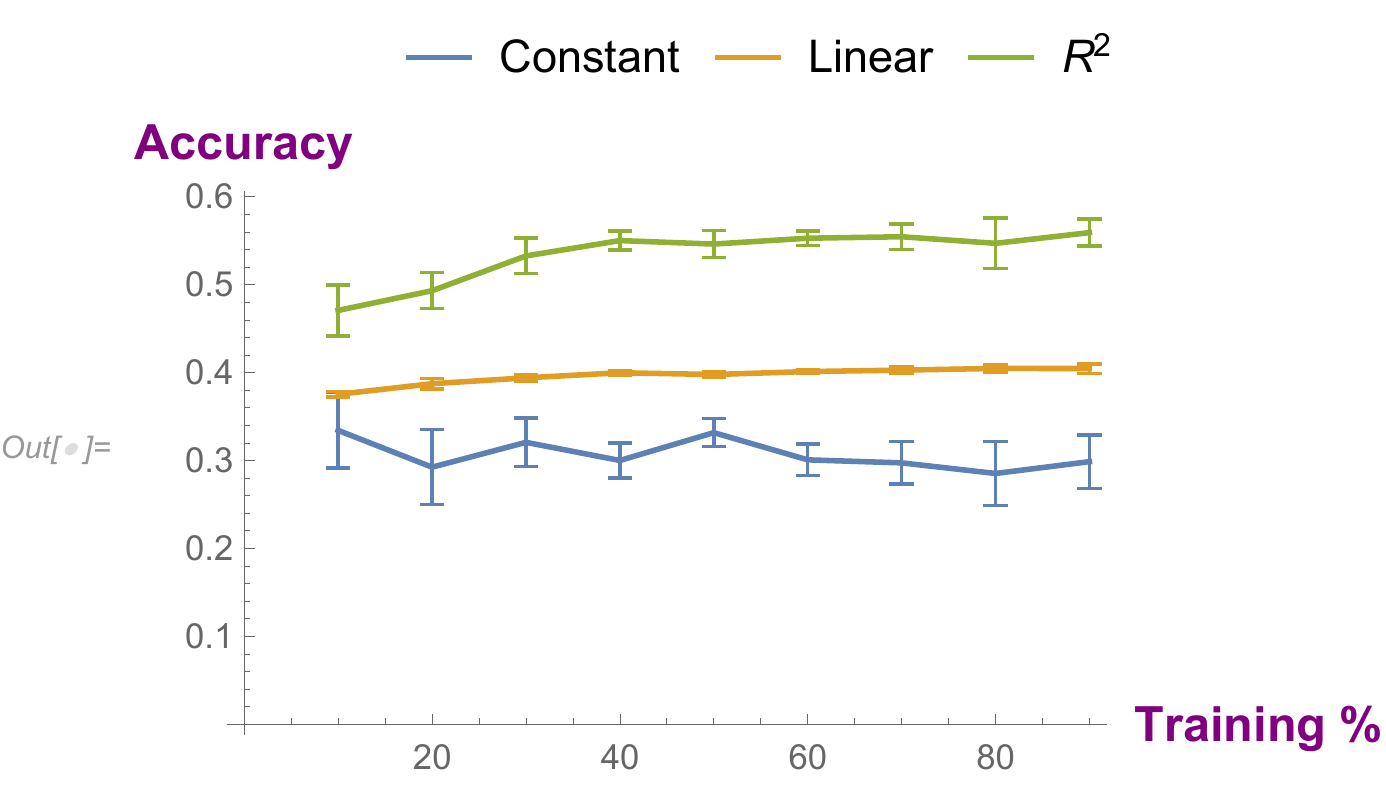}
}
\caption{{\sf {\small
The training curve for the neural network in \eqref{NN} for ML of the maximum Laplacian eigenvalue of our graph dataset.
A random sample of $x\%$ is trained and then the NN is validated on the remaining $(100-x)\%$; we let x range from 10 to 90 in increments of 10.
We show 3 measures of goodness of fit: the constant and linear terms of the linear regression of predicted and actual values, as well as the coefficient of determination $R^2$.
}}
\label{f:minDeltaEigenTC}}
\end{figure}

We use a slightly improved NN than \eqref{NN} by inserting an element-wise hyperbolic tangent layer between the linear layer and the final output summation layer.
Still, the performance is not as good as that of $L$ is \S\ref{s:gap}.
We present the linear fit between the actual and the predicted value on the validation set in part (a) of Figure \ref{f:minDeltaEigenTC}.
The line is found to be
\begin{equation}
y \simeq 0.884 x + 0.049 \ , \qquad R^2 \simeq 0.753 \ .
\end{equation}
 In part (b) of the figure, we present the training curve from 10\% to 90\% training data as before and plot the behaviour of the linear and constant terms of the linear regression fit, as well as the $R^2$-coefficient.

As a double check against \S\ref{s:maxL}, we also performed the above supervized ML for the maximal eigenvalue of normalized Laplacian $\Delta$. Again, the performance, though satisfactory, is not as impressive as \eqref{fit20maxEigen}.
Here, at 20-80\% split, we have $y \simeq 0.816 x + 0.314$ with $R^2 \simeq 0.575$.

It is curious that eigenvalues of the ordinary Laplacian $L$ seem more amenable to ML than our normalized Laplacian $\Delta$; one reason for this that the former has much larger variation whilst the latter is confined to within $[0,2]$ and our NN structure is doing non-linear regression without attempting to detect scale.

~\\
~\\

\section{Ricci Flat Graphs}\setall
As with Riemannian manifolds, once a notion of curvature is established, one of the first explorations to undertake is Ricci-flatness.
For K\"ahler manifolds, this led to the Calabi Conjecture/Yau Theorem.
For (locally finite) graphs, continuing with the idea of $\Gamma$ in Definition \ref{defRicci}, a curvature-dimension type inequality was introduced \cite{ly}:
\begin{definition}\label{CDkm}
A curvature-dimension inequality can be introduced to the graph Laplacian $\Delta$ as
\[
\Gamma_2(f,f)(x) \geq \frac{1}{m} \big( \Delta f (x) \big)^2 + k(x) \Gamma(f,f) (x) \ ,
\]
where $m \in \IR_{\geq 1}$ is called the {\em dimension} of $\Delta$ as an operator and $k(x)$ is the lower bound on the Ricci curvature associated to $\Delta$.
Let us call this inequality type $CD(m,k)$.
\end{definition}

The next concepts we need are distributions and transport \footnote{
The transport distance is usually first defined in terms of a coupling but we will adhere to the one thus presented to minimize introducing new concepts.
} on a (locally finite) graph $G=(V,E)$:
\begin{definition}
A {\bf probability distribution} over $V$ is a map $\mu : V \to [0,1]$ such that $\sum\limits_{x \in V} \mu(x) = 1$.
A function $f : V \to \IR$ is called {\bf 1-Lipschits} if $f(x) - f(y) \leq d(x,y)$ for each $x,y \in V$, i.e., the difference in value of the function between any two vertices is bounded by the distance between them.
The {\bf transportation distance} between two probability distributions $\mu_1$ and $\mu_2$ is
\[
W(\mu_1, \mu_2) := \sup \limits_{f \mbox{ {\tiny 1-Lipschits}}} \sum\limits_{x \in V} f(x) \big( \mu_1(x) - \mu_2(x) \big) \ .
\]
\end{definition}

\subsection{Two Notions of Ricci-Flatness for Graphs}
Thus prepared, we have two notions of Ricci-flatness for graphs and we will proceed with care to compare them
\footnote{There is actually a third notion of Ricci-flatness, which we will not need in this paper.
This definition requires the concept of a $k$-frame as was introduced by \cite{fy}; it also requires that the graph be {\bf regular} in that all vertices have the {\it same degree/valency}.
We say that a regular graph $G = (V,E)$ has a {\bf  local $k$-frame} at vertex $x$ if there exist injective mappings $\eta_{i = 1, \ldots k}$ from a neighbourhood of $x$ into V such that (1) $x$ is adjacent to $\eta_i x$ for all $1 \leq i \leq k$ and (2) $\eta_i x \neq \eta_j x$ for $i \neq j$.
Then, we can define the regular graph $G$ to be Ricci flat at $x$ if there is a local $k$-frame in the neighbourhood of $x$ such that for all $i = 1, \ldots, k$ we have a notion of commutativity of transport:
\[
\bigcup\limits_{j=1}^k (\eta_j \eta_i) x =
\bigcup\limits_{j=1}^k (\eta_i \eta_j) x \ .
\]
One might refer to this third notion as $k$-frame-Ricci-flatness.
}
.
In some sense, Ricci-flat graphs are the discrete analogues of Calabi-Yau manifolds.
It is still an open conjecture by the second author that in each complex dimension $n$, there is a finite number (in the sense of  topological type) of compact, smooth Calabi-Yau $n$-folds \footnote{In complex dimension $n=1$, for instance, there is only $T^2$. For $n = 2$, there is only $T^3$ and the K3 surface. For $n=3$, the discovered number is huge (on the order of $10^{10}$) but the number of topologically inequivalent CY 3-folds is expected to be finite; cf.~brief review in \cite{Bao:2020sqg}.}

One way to define Ricci-flatness is to modify the notion of Ricci curvature for metric spaces in the sense of \cite{oll} 
The classification was addressed in \cite{lly} (cf.~also \cite{ckllly,ckllly2,osy}).
For vertices $x,y \in V$, and any real value $\alpha \in [0,1]$, define the probability distribution $\mu_x^\alpha$ and using the transportation distance $W$:
\begin{equation}
\mu_x^\alpha(z) := \left\{
\begin{array}{ccl}
\alpha \, && \mbox{ if } z = x \ , \\
\frac{1 - \alpha}{d(x)} \ , && \mbox{ if } z \sim x \ , \\
0 && \mbox{ otherwise \ ; }
\end{array}
\right. \qquad
k_\alpha(x,y) := 1 - \frac{W(\mu_x^\alpha, \mu_x^\alpha)}{d(x,y)} \ .
\end{equation}
Then, we have an Ollivier-type Ricci curvature as defined in \cite{lly}:
\begin{definition}
Define curvature $k(x,y) = \lim\limits_{\alpha \to 1} \frac{k_\alpha(x,y)}{1 - \alpha}$.
A (locally finite) graph is Ricci-flat if $k(x,y)$ vanishes for any edges $x \sim y$.
\end{definition}
We will refer to this notion as {\bf OLLY-Ricci-flatness}.

Another way to define graph Ricci-flatness is to consider the inequality $CD(m,k)$ in Definition \ref{CDkm} whose classification was addressed in \cite{HuaLin}:
\begin{definition}
A Ricci-flat graph $G$  in the sense of curvature-dimension is one whose Laplacian $\Delta$ satisfies $CD(\infty, 0)$.
\end{definition}
We will refer to this definition as {\bf CD-Ricci-flatness}.

\subsubsection{Classification Results}
As mentioned, both notions have undergone classification.
First, a key result of \cite{lly} is that whilst for girth (cf.~definition in \eqref{girth}) 3 and 4, there are an infinite number of OLLY-Ricci-flat graphs (q.v.~\cite{bly}), for girth 5 or more, we have that \footnote{Note that the initial classification in \cite{lly} missed the Triplex graph, which we later found in \cite{ckllly2}.}
\begin{theorem}[Lin-Lu-Yau]
If $G$ is a (locally finite) OLLY-Ricci-flat graph with girth $\geq 5$, then it is one of the 2 infinity families:
\begin{itemize}
\item the infinite line;
\item the cycle graph $C_{n \geq 6}$;
\end{itemize}
or one of the 4 exceptional cases:
\[
\includegraphics[trim=0mm 0mm 0mm 0mm, clip, width=5in]{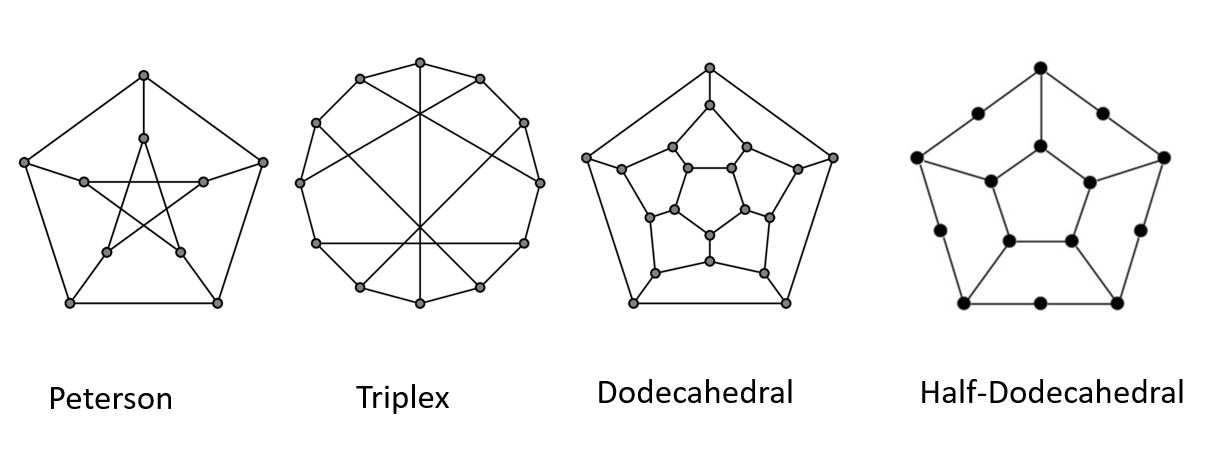}
\]
\end{theorem}
It is further curious that the classification is reminiscent of an ADE pattern: a simple infinite family, a more complicated infinite family, and coincidentally 3 exceptional cases (if we combined the dodecahedral with the half-dodecahedral).

Second, the main result of \cite{HuaLin} is the classification of CD-Ricci-flat graphs (in cit.~ibid., the authors consider weighting and normalization so their results are more general.
For our present purposes, the vertex weights $m$ and the edge weights $\mu$ are all set to 1 and their theorem 1.3 applies):
\begin{theorem}[Hua-Lin]\label{thm:CDricci}
Let $G=(V,E)$ be of girth at least 5, then it satisfies $CD(0,\infty)$ and is thus CD-Ricci-flat if it is one of the following:
\begin{itemize}
\item The path graph $P_{k \geq 1}$, the cycle graph $C_{n \geq 5}$;
\item The infinite line $P_{\IZ}$ or the infinite half-line $P_{\mathbb{N}}$;
\item The star graphs $\mbox{Star}_{n \geq 3}$ 
\[
\begin{array}{c}
\includegraphics[trim=0mm 0mm 0mm 0mm, clip, width=0.8in]{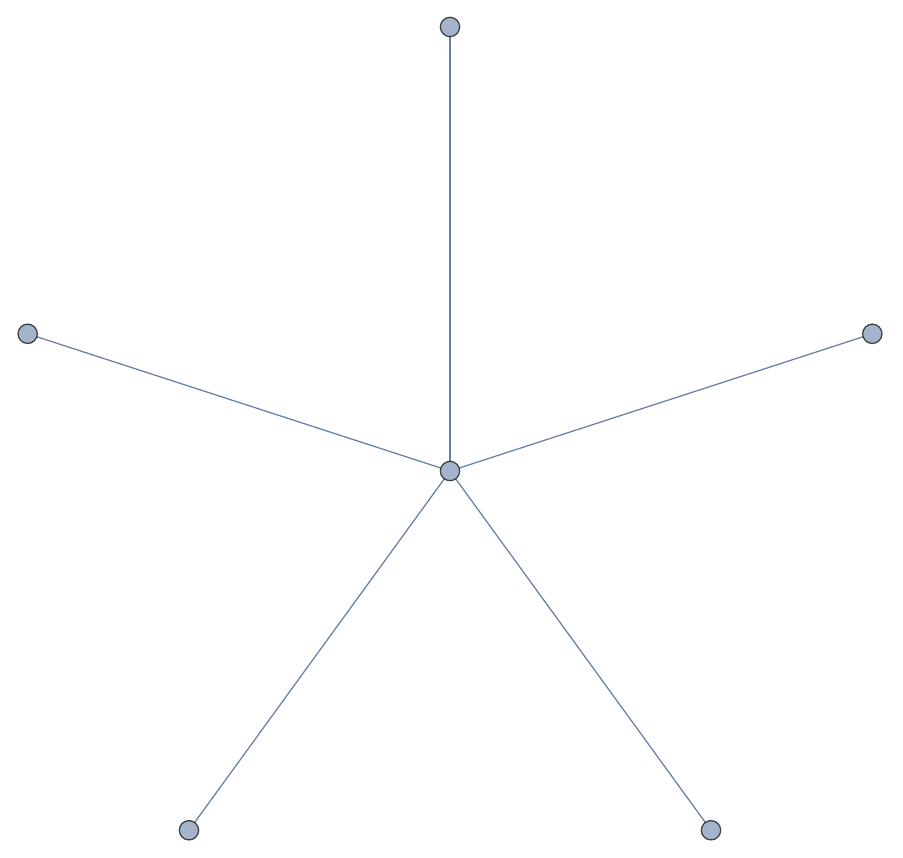}
\mbox{Star}_5 \ ,
\includegraphics[trim=0mm 0mm 0mm 0mm, clip, width=0.8in]{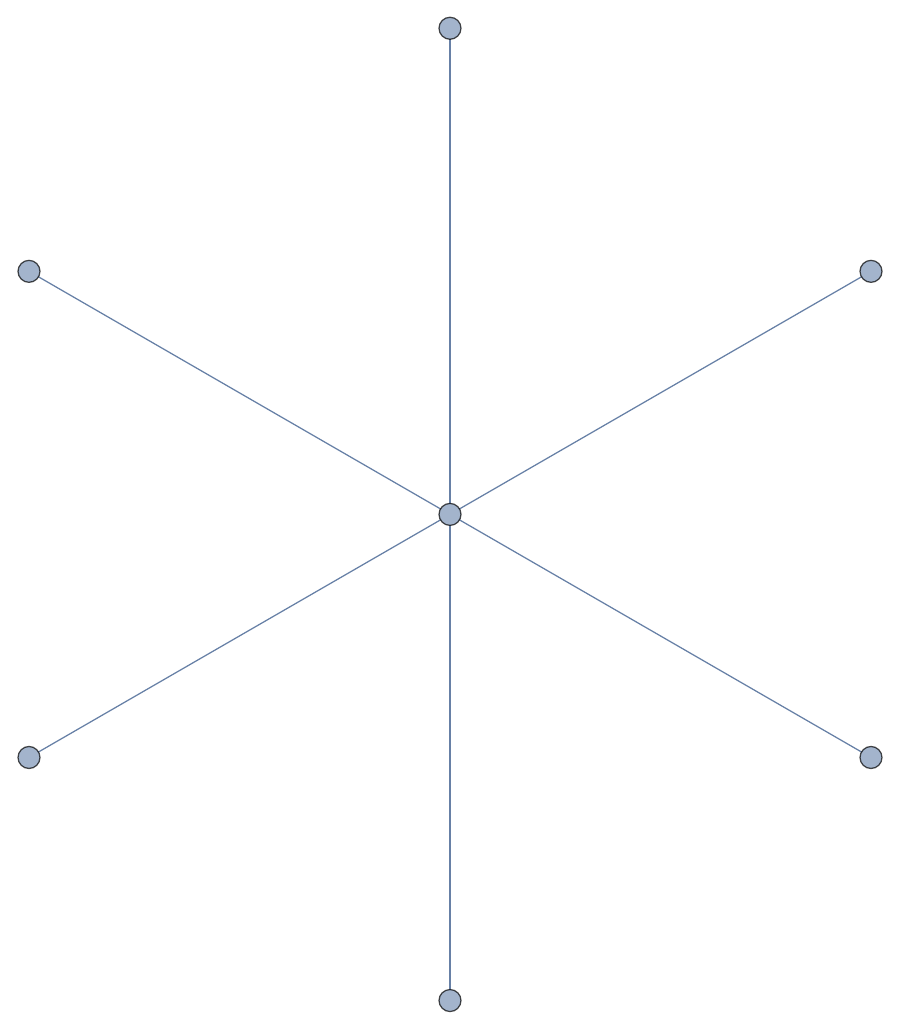}
\mbox{Star}_6 \ ,
\includegraphics[trim=0mm 0mm 0mm 0mm, clip, width=0.8in]{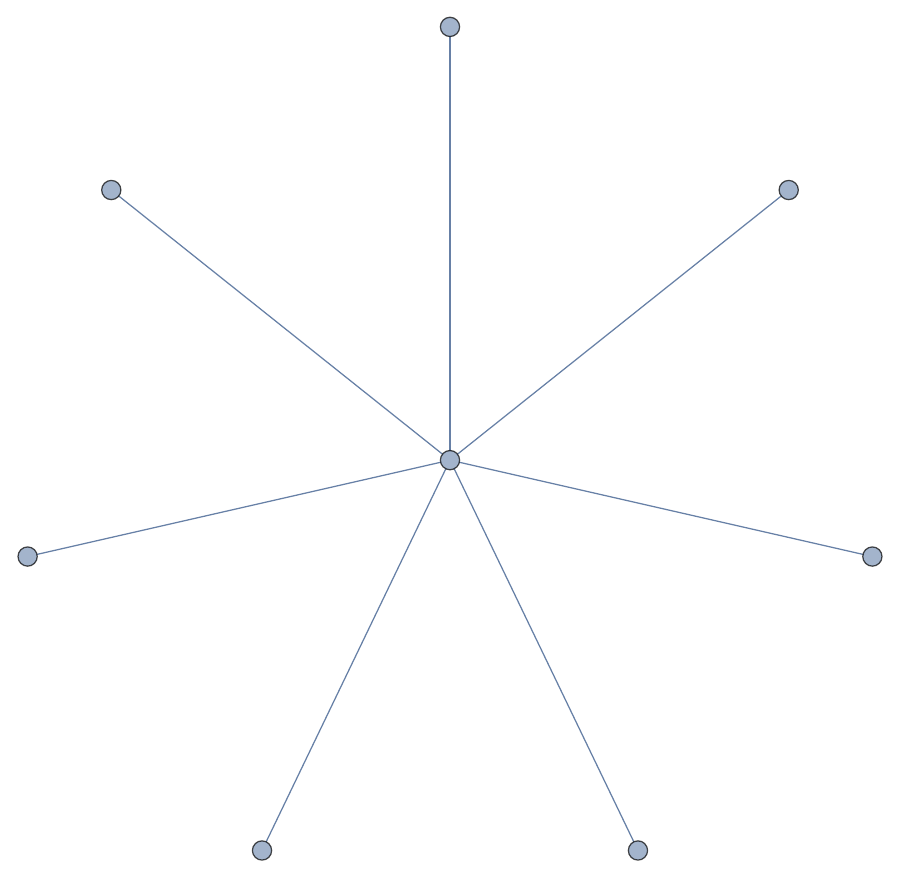} 
\mbox{Star}_7 \ , 
\ldots \ldots
\end{array}
\]
\item The extended star graphs $\mbox{Star}_{3}^{i=1,2,3}$:
\[
\begin{array}{c}
\includegraphics[trim=0mm 0mm 0mm 0mm, clip, width=4in]{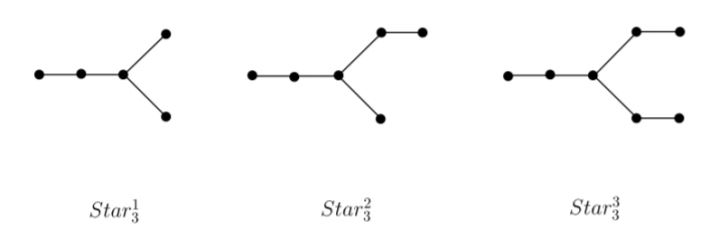}
\end{array}
\]
\end{itemize}
\end{theorem}
Thus one again, we have a somewhat ADE-type pattern with an infinite family of simple cases (together with the 1/2-infinite versions), an infinite family of less-simple cases (the star graphs), as well as 3 exceptionals (incidentally, the three extended star graphs are precisely the Dynkin diagrams for $D_5$, $E_6$ and $\widehat{E_6}$).

\paragraph{Question: }
A question immediately comes to mind as to whether AI can distinguish a Ricci flat graph just by looking at it.
A similar venture was undertaken in \cite{He:2017aed,He:2019vsj,He:2019nzx}.
In particular, in \cite{He:2019vsj}, the question of whether ML can tell if a Calabi-Yau manifold is elliptically fibred was posed, and answered in the affirmative.
There, by going over the data-set of complete intersection Calabi-Yau manifolds in products of projective spaces, wherein the elliptic fibrations have been found using tradition techniques in algebraic geometry \cite{Anderson:2017aux}, a NN was set up and was found to over 99\% confidence that such fibration structures can be machine-learned.

Now, non-elliptic manifolds are actually quite rare in the space of Calabi-Yau manifolds and data enhancement in the manner of what we will shortly describe was first performed.
The situation is similar here:
Ricci-flat graphs in both senses of the definitions are relatively rare in the space of graphs, especially in the OLLY sense.
We will thus enhance the data by performing appropriate permutations of the adjacency matrices, which are clearly equivalent representations.

\begin{figure}[!ht!b]
{
(a)
\includegraphics[trim=12mm 0mm 0mm 0mm, clip, width=2.5in]{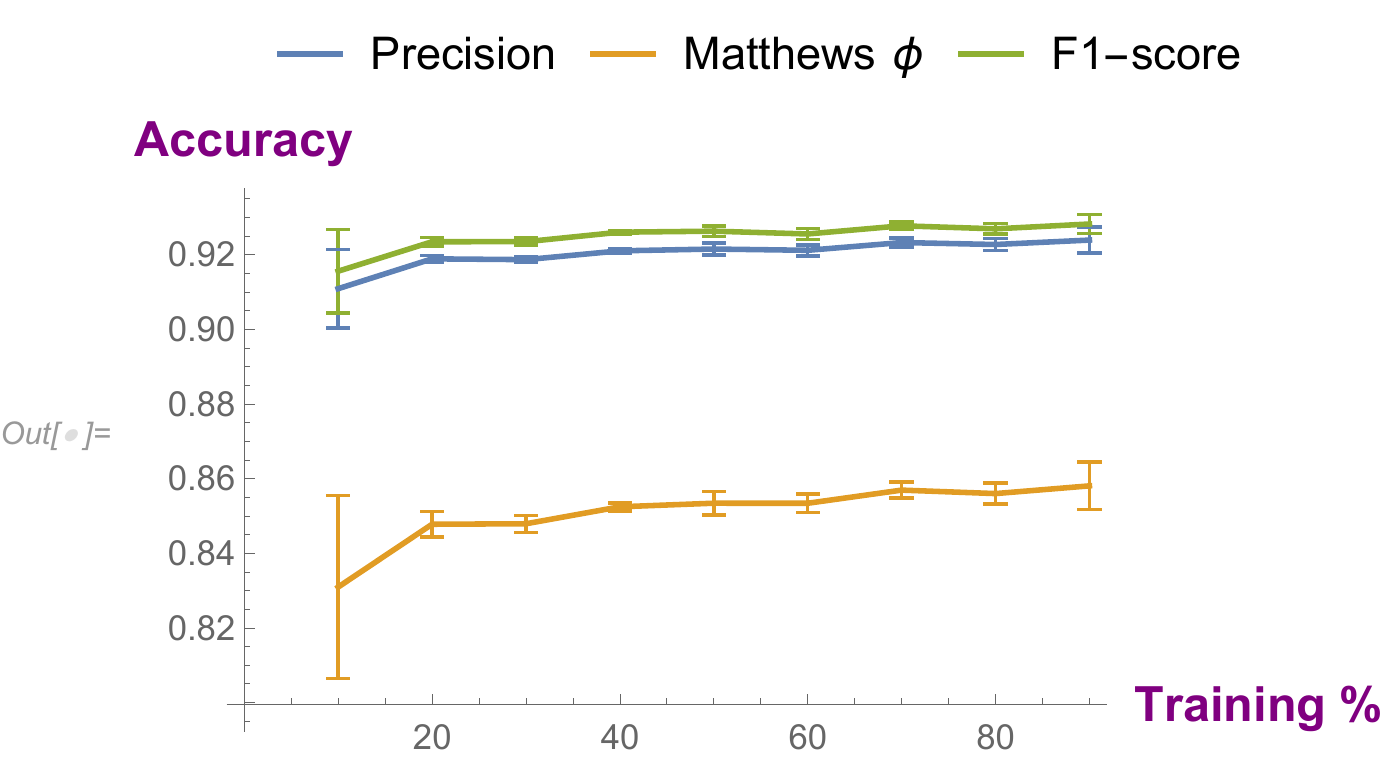}
(b)
\includegraphics[trim=12mm 0mm 0mm 0mm, clip, width=2.5in]{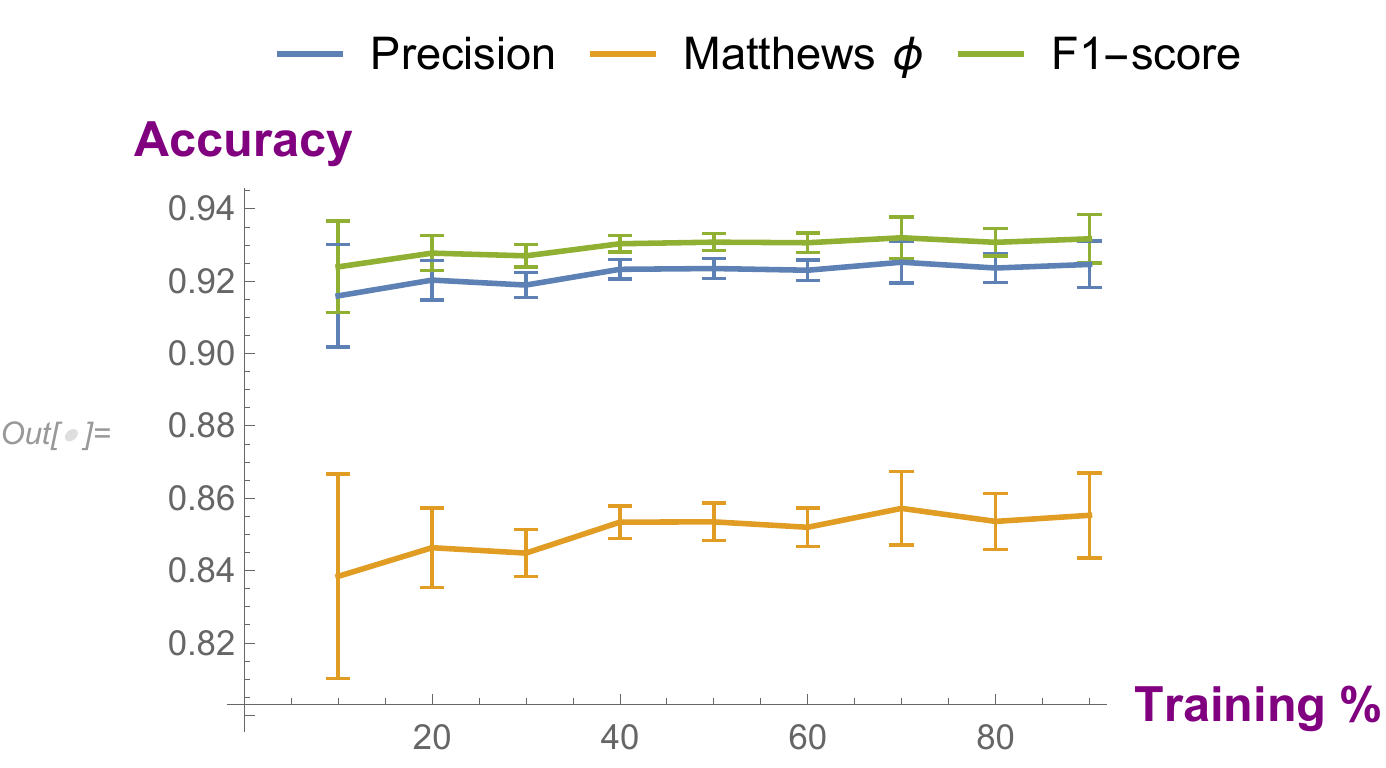}
}
\caption{{\sf {\small
The training curve for a classifier (optimized between a logistic regression and a gradient boosted decision tree) for distinguishing a Ricci-flat graph from our dataset, respectively for the OLLY and CD notions of the RIcci-flatness in parts (a) and (b).
A random sample of $x\%$ is trained and then the NN is validated on the remaining $(100-x)\%$; we let x range from 10 to 90 in increments of 10.
We show 3 measures of goodness of fit: the naive precision, the F1-score and the Matthews coefficient.
}}
\label{f:RicciFlat}}
\end{figure}

\subsection{Distinguishing OLLY-Ricci-Flat Graphs with ML}
In our dataset, we have 1806 graphs of girth at least five \footnote{There are 1805 built in, including the cycle graphs, the dodecahedron, the Triplex, and the Petersen, but the half-dodecahedron is added in by hand.}.
Amongst these let us consider the OLLY-Ricci-flat ones: the infinite line is out since we only consider finite graphs, so we have the cycles graphs and the 4 exceptional cases, totaling 20.

We enhance the data by permuting the adjacency matrices of the Ricci-flat cases by 800 times randomly, and that of the non-Ricci-flat cases by 8 times randomly.
This gives us a rather balanced labeled set of around 15000 each of ``1'' (Ricci flat) and ``0" (not Ricci flat), perfectly adapted for a binary classification.

At 5-fold cross-validation and using a optimzed classifier of gradient-boosted decision trees and logistic regression, we find that
\begin{equation}
P_N \simeq 0.925 \pm 0.004, \ F_1  \simeq 0.927 \pm 0.004, \ \phi \simeq  0.855 \pm 0.007\ ,
\end{equation}
which is excellent.
The training curve is shown in part (a) of Figure \ref{f:RicciFlat}.
We see that distinguishing OLLY-Ricci-flat graphs is performed very well by ML, in that we are into the 90's in terms of measures of goodness of fit.
As an extra precaution, what if we {\it randomly} assigned 1 to around 50\% of the data and 0 to the remainder?
The 5-fold cross validation gave
$P_N \simeq 0.546\pm0.029$, $F_1 \simeq 0.101 \pm 0.205$ and $\phi \simeq 0.004 \pm 0.0135$.
This means that the precision is at around $50\%$, as good as random guessing, which is further confirmed by the $F_1$-score and Matthews $\phi$-coefficient being near 0.
This is re-assuring indeed and shows that there is a significant {\it pattern} which distinguishes Ricci-flatness, and that the pattern is not random.

\subsection{Distinguishing CD-Ricci-Flat Graphs with ML}
Having gained confidence, let us move onto the other notion of Ricci-flatness.
Within our dataset of 1805 graphs with girth at least 5, there are 58 which satisfy the criterion of CD-Ricci-flatness in accordance with Theorem \ref{thm:CDricci}.
We thus perform 300 random permutations of the adjacency matrices for the ``yes'' cases and 5 for the ``no'' cases, giving us about 9000 each, of a roughly balanced 18000 database ready for binary classification.

At 5-fold cross-validation and using a optimzed classifier of gradient-boosted decision trees and logistic regression, we find that
\begin{equation}
P_N \simeq 0.924 \pm 0.004, \ F_1  \simeq 0.931 \pm 0.004, \ \phi \simeq  0.854 \pm 0.006\ ,
\end{equation}
which is again excellent.
The training curve is shown in part (b) of Figure \ref{f:RicciFlat}.

We conclude that distinguishing Ricci-flat graphs, in both the OLLY and the curvature-dimension notions, is performed comparably well by ML, in that we are into the 90's in terms of measures of goodness of fit by an optimized binary classifier.
The high percentage agreement is further checked by the high F1-score and Matthews coefficient, showing that the classification is truly good and that false-positives and false-negatives are insignificant.

\section{Homology and Cohomology of Graphs}\setall
A fundamental result in Riemannian geometry is the decomposition of Hodge which implies that (co-)homology of the manifold should be described by the zero-modes of the Laplacian.
Can this be carried over to our graphical context?
In \cite{glmy1,glmy2}, a notion of differential forms and boundary operators, and hence cohomology and homology, were introduced on {\bf directed} simple graphs.

In brief, the construction proceeds as follows.
\begin{itemize}
\item Let $G = (V,E)$ be a finite directed graph (digraph) so that $V = \{i\}_{i = 1, \ldots, n}$ are the vertices 
and $E = \{ i_k i_{k+1} \}$ are directed edges (arrows) from vertex $i_k$ to $i_{k+1}$.
On $G$, we have
	\begin{itemize}
	\item {\bf Elementary $p$-path} $e_{i_0 \ldots i_p}$, which is  any ordered sequence of $p+1$ vertices $i_0$, $i_1, \ldots$ to $i_p$;
	\item Elementary {\bf regular} $p$-path is one for which $i_k \neq i_{k+1}$, so that there is no back-tracking; 
	\item {\bf Allowed} elementary $p$_path is one for which $i_k \to i_{k+1} \in E$ for all $k = 0, \ldots, p-1$, i.e., the path is actually traversed by arrows in $G$; 
	\end{itemize}

\item Fix $\IK$ to be a commutative ring with unity (we mostly just take $\IR$) and consider the free $\IK$-module generated by the elementary $p$_paths:
	\begin{equation}
	\Lambda_p := \{ \sum k e_{i_0 \ldots i_p} \} 
	=
	\mbox{Span}_{\IK} \{  \mbox{elementary $p$-paths} \}
	\ ;
	\end{equation}
	elements of $\Lambda_p$ are just called $p$-paths.
		\begin{itemize}
		\item Define the boundary operator $\partial: \Lambda_{p+1} \to \Lambda_p$ by
		\begin{equation}
		\partial e_{i_0 \ldots i_p} := \left\{
		\begin{array}{l}
		\sum\limits_{q=0}^p (-1)^q e_{i_0 \ldots \hat{i_q} \ldots i_p} \mbox{ for } p \geq 1 \ , \\
		0 \mbox{ for } p = 0 \ ,
		\end{array} \right.
		\end{equation}
		where the hat means omission and we set $\Lambda_{-1} = \{0\}$;
		\item Then we have $\partial^2 = 0$;
		
		\item Define the subspace of regular paths
		\begin{equation}
		\cR_p := \mbox{Span}_{\IK} \{  \mbox{regular elementary $p$-paths} \} \subset \Lambda_p \ ,
		\end{equation}
		which is linearly isomorphic to the quotient of $\Lambda_p$ by the irregular paths; i.e., 
		we set all irregular paths to 0.
		
		\item Define one further step the space of allowed paths
		\begin{equation}
		\cA_p := \mbox{Span}_{\IK} \{  \mbox{allowed regular elementary $p$-paths} \} 
			\subset \cR_p \subset \Lambda_p
		\end{equation}
		and consider the $\partial$-invariant subspaces (i.e., the boundaries are still allowed paths)
		\begin{equation}
		\Omega_p := \{ v \in \cA_p : \partial v \in \cA_{p-1} \} \subset \cA_p \ .
		\end{equation}
		By construction, the 0 and 1-paths are just the vertices and arrows respectively:
		\begin{equation}
		\Omega_0 = \cA_0 = V, \quad \Omega_1 = \cA_1 = E \ ,
		\end{equation}
		so $\dim \Omega_0 = |V|$ and $\dim \Omega_1 = |E|$.
		Then, we have the chain complex
		\begin{equation}
		0 \leftarrow \Omega_0 \stackrel{\partial}{\leftarrow} \Omega_1 \stackrel{\partial}{\leftarrow} \ldots
		\stackrel{\partial}{\leftarrow} \Omega_p \stackrel{\partial}{\leftarrow} \ldots 
		\end{equation}
		from which can define {\bf graph homology groups} $H_p(G) := \ker(\partial_p) / \im(\partial_{p+1})$.
		\end{itemize}

\item Next, we can define the dual to the above. First, a $p$-form on $G$ is a $\IK$_valued function on $V^{p+1}$, i.e., it is a function $\omega(i_0, \ldots, i_p)$ with $p+1$ arguments.
We have the freely generated $\IK$-module
\begin{equation}
\Lambda^p :=  \mbox{Span}_{\IK} \{ \mbox{$\IK$-valued function $\omega$} \} =
\{
\sum\limits_{i_0, \ldots, i_p \in V} \omega(i_0, \ldots, i_p) e^{i_0 \ldots i_p)}
\} \ ,
\end{equation}
where $e^{i_0 \ldots i_p)}$ is a canonical basis of $\Lambda^p$.
	\begin{itemize}
	\item On $\Lambda^p$ we can define an exterior derivative $d : \Lambda^p \to \Lambda^{p+1}$ as
	\begin{equation}
	(d \omega) (i_0, \ldots i_{p+1}) := \sum\limits_{q=0}^{p+1} (-1)^q \omega(i_0, \ldots, \hat{i_q}, \ldots i_p) \ ,
	\end{equation}
	where the hat is again omission of the index.
	We can check that $d^2 = 0$.
	
	\item Again, regularity can be defined by having no back-tracking of indices: $i_k \neq i_{k+1}$, so that we have the subspace of regular $p$-forms
	\begin{equation}
	\cR^p :=  \mbox{Span}_{\IK} \{  e^{i_0 \ldots i_p} : i_k \neq i_{k+1} \} \subset \Lambda^p \ .
	\end{equation}
	
	\item We now wish to quotient out by the non-allowed $p$_forms, corresponding to the non-allowed $p$-path:
	\begin{equation}
	\cN^p := \mbox{Span}_{\IK} \{  e^{i_0 \ldots i_p} :  i_0 \ldots i_p \in \Lambda_p \backslash \cA_p \} \ .
	\end{equation}
	
	\item Subsequently we can define $\Omega^p := \cR^p / (\cN^p + d \cN^{p-1})$ by treating non-allowed $p$-forms as zero and check that $\Omega^p$ is $d$-invariant. This gives us the dual complex
	\begin{equation}
	0 \rightarrow \Omega^0 \stackrel{d}{\rightarrow} \Omega^1 \stackrel{d}{\rightarrow} \ldots
		\stackrel{d}{\rightarrow} \Omega^p \stackrel{d}{\rightarrow} \ldots 
	\end{equation}
	from which we can define {\bf graph cohomology groups} $H^p(G) = \ker(d_{p}) / \im(d_{p-1})$.
	\end{itemize}
\end{itemize}

\comment{Happy new year ! 
I wrote Some papers on defining cohomology for direct graphs 
They are interesting 
We also did works in Ricci flat graphs
Yau
}

As usual with duality between homology and cohomology (of compact smooth manifolds), a non-degenerate bilinear pairing can be established between $H_p(G)$ and $H^p(G)$, rendering them isomorphic as dual vector spaces.
In particular, they have the same dimension.
Note that
\begin{align}
\nn
&\dim H_p(G) = \dim \Omega_p - \dim \partial \Omega_p - \dim \partial \Omega_{p+1} = \\
\label{dimH}
&\dim H^p(G) = \dim \Omega^p - \dim d \Omega^p - \dim d \Omega^{p-1} \ . 
\end{align}
Whence a graph Euler number can be defined as
\begin{equation}
\chi(G) = \sum\limits_{p=0}^n (-1)^p \dim H_p(G) \ ,
\end{equation}
with $n$ sufficiently large so that $\dim H_{p > n}(G) = 0$.
Due to the alternating sum in \eqref{dimH}, as always, we have an Euler-Poincar\'e type of relation
\begin{equation}
\chi(G) = \sum\limits_{p=0}^n (-1)^p \dim \Omega_p(G) \ . 
\end{equation}

\subsection{Machine Learning Graph Euler Number}
We leave a detailed calculation of graph homology to the Appendix, demonstrating its non-triviality.
For manifolds, the application of machine-learning to topological invariants was initiated in \cite{He:2017aed}.
It was found, for example, that Hodge numbers of classes of Calabi-Yau maifolds respond well to a simple neural network.
Naturally, one could ask whether a similar behaviour occurs here.

First, we need to establish a reasonable dataset.
We can take our database of undirected graphs, and assign random directions since the above notion of homology requires orientation.
To be precise, we take a selection of 10 undirected graphs (up to 6 vertices since the homology computation gets quite intensive) in the set, and assign 100 random directions to each, and randomly permute the $\chi=0$ (which are slightly under-represented) adjacency matrices 48 times and the $\chi \neq 0$ cases 20 times.
This gives us a binary classification problem of around 15K in each category.
Interestingly, neither a classifier (optimized between logistic regression and a decision tree) nor a neural network of the type in \eqref{NN} could do this problem well: we obtain naive precision around 0.50 with $\phi \simeq 0.08$, which is marginally better than random guessing.
We tried other NNs such gated recurrence networks, which performed slightly better at $\phi \simeq 0.15$ but could not find any ML which did so nicely as distinguishing properties such as Ricci-flatness, planarity, chromaticism, etc, as above.

\begin{figure}[!ht!b]
\centerline{
(a)
\includegraphics[trim=12mm 0mm 0mm 0mm, clip, width=3in]{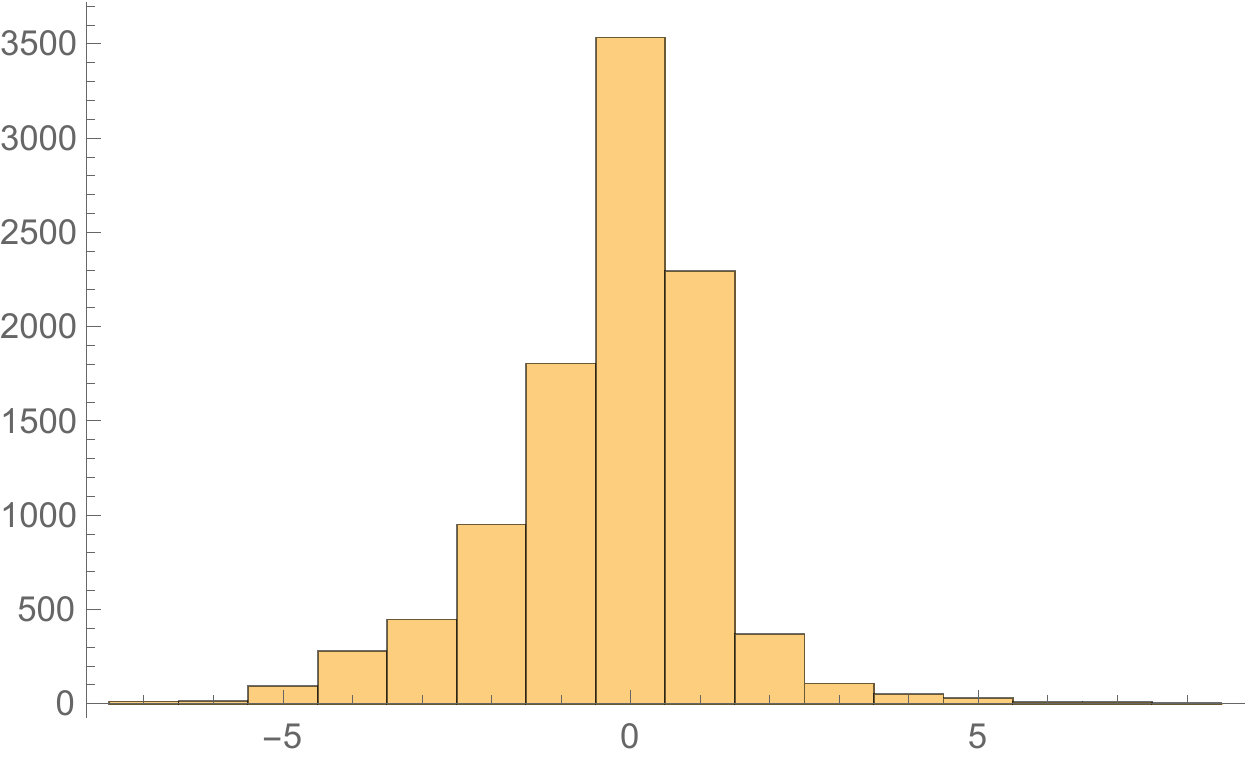}
(b)
\includegraphics[trim=12mm 0mm 0mm 0mm, clip, width=2in]{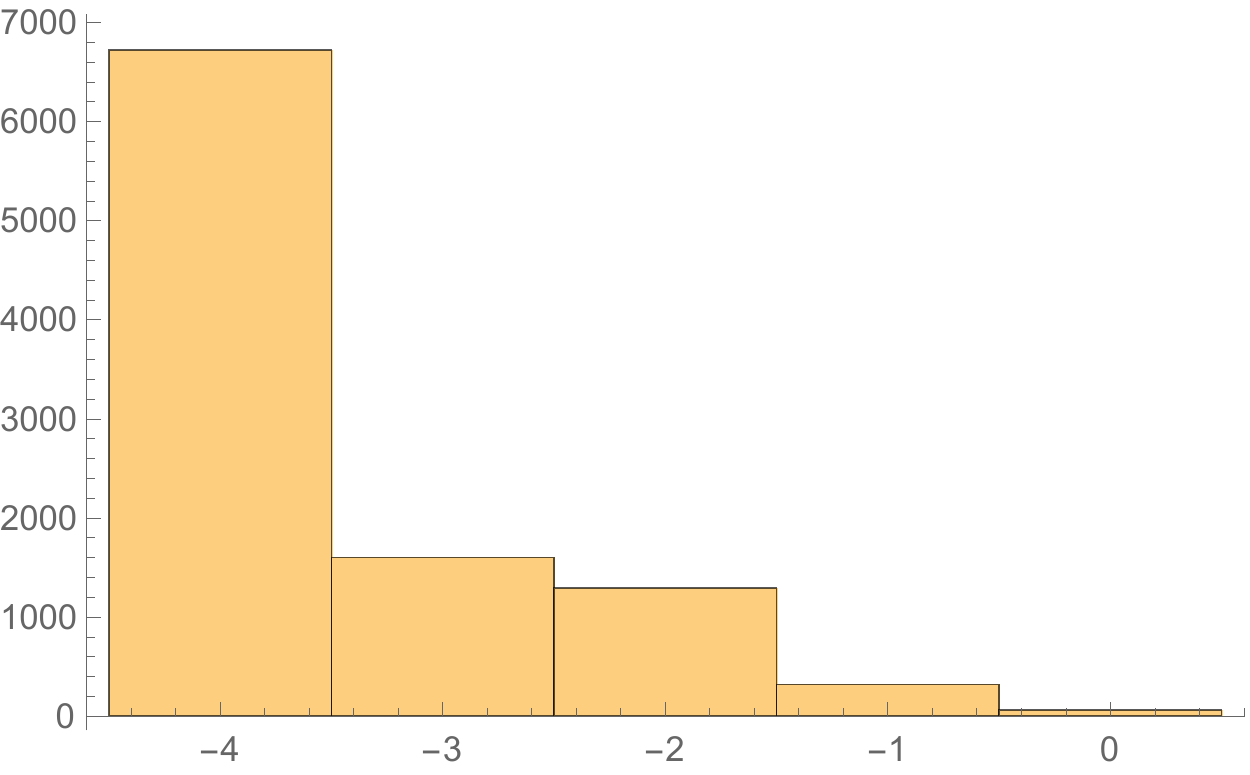}
}
\caption{{\sf {\small
(1) histogram for $\chi$ of 10K random graphs with random orientations;
(2) histogram for $\chi$ of 10K random orientations on the cube. 
}}
\label{f:euler}}
\end{figure}

It is curious that learning such topological quantities for graph did not behave as brilliant as the counterpart in manifolds \cite{He:2017aed,Bull:2018uow}.
For reference, we give an idea of the distribution of graph Euler number in Figure \ref{f:euler}.
In part (a), we take the first 100 graphs in the dataset, with 100 random direction assignments to each, giving us 10K graphs for which we compute the Euler number and show the distribution, a somewhat Gaussian around 0.
In part (b), we take a particular graph with 8 vertices, say the cube, and then assign 5K random orientations to compute the Euler numbers, and show the distribution.

\section{Conclusions and Prospectus}
Due to their combinatorial nature and hence the facility of representation as matrix manipulations, finite graphs present a natural venue for a diverse number of fields in mathematics and physics.
Inspired by a long programme of the second author to investigate discrete analogues of the foundational problems in Riemannian geometry via finite graphs, as well as a more recent programme of the first author to explore how different mathematical structures can be machine-learned or detected by generic neural networks without a priori knowledge, we have employed some of the latest techniques in data science to study aspects of graph theory in connection to Riemannian geometry.

Taking the freely available Wolfram database of finite simple graphs up to 100 nodes as a concrete playground - a representative set of some 8000 graphs - we have examined a host of relevant statistical and machine-learning properties.
We warmed up with machine-learning a miscellany of basic graph quantities in a supervised learning paradigm where we trained a classifier/neural network with the labeling ``adjacency matrix $\longrightarrow$ property''.
These are mostly done with random-forest decision trees and logistic regressions that have no knowledge of the underlying graph theory.
Denoting the triple-check for accuracy measures as (naive precision, F1-score, Matthews $\phi$ coefficient) in the case of binary classification or a double-check as (naive precision, Matthews $\phi$ coefficient) in the case of multi-category classification, we summarize the 5-fold cross-validation results in Table \ref{t:summary}.
For clarity of comparison we order from best performance downward.

Next, we attempted a more refined predictor to study the maximal and minimal eigenvalues of the Laplacian.
Using a simple neural network with the architecture given in \eqref{NN}, with standard activation functions, 
we attempted to predict these eigenvalue bounds by looking solely at the adjacency matrix.
The performance is again very good, in that the actual and predicted values fit to a line $y = 1.01 x - 0.06$ with $R^2 = 0.93$ (perfect prediction would mean $y = x$ with $R^2= 1$) for the maximal eigenvalue.
Interestingly, the minimal positive eigenvalue (the spectral gap) behaved a bit worse in the prediction, with $y = 0.91 x + 0.70$ with $R^2 = 0.897$.

\begin{table}[t!!!]
\begin{tabular}{c|c}
Property & Accuracy Measure \\ \hline
whether graph is acyclic & $(0.954 \pm 0.001 \ , 0.955 \pm 0.001 \ , 0.912 \pm 0.002)$ \\
whether graph is OLLY-Ricci flat & $(0.925 \pm 0.004 \ , 0.927 \pm 0.004, \ , 0.855 \pm 0.007)$\\
whether graph is CD-Ricci-flat & $(0.924 \pm 0.004 \ , 0.931 \pm 0.004 \ , 0.854 \pm 0.006)$ \\
whether graph is genus $<=0$, $=0$ or $>0$ & $(0.814 \pm 0.003 \ , 0.721 \pm 0.005)$ \\ 
whether graph is planar & $(0.812 \pm 0.004, \ 0.832 \pm 0.004, \ 0.619 \pm 0.009)$ \\
girth of graph $=2$, $=4$ or $>4$ & $(0.771 \pm 0.017  \ , 0.656 \pm 0.026)$ \\
diameter $\leq 2$, $=3,4$ or $> 4$ & $(0.765 \pm 0.004 \ ,  0.647 \pm 0.005)$ \\ 
skewness of graph $=0$, $=1$ or $>1$ & $(0.747 \pm 0.005$ and $0.597 \pm 0.008)$ \\
whether there is a Hamilton cycle & $(0.781 \pm 0.008 \ , 0.770 \pm 0.009 \ , 0.564 \pm 0.017)$\\
whether there is a Eulerian cycle & $(0.731 \pm 0.015\ , 0.721 \pm 0.026 \ , 0.473 \pm 0.024)$\\
Round of Laplacian eigenvalue & $(0.4982 \pm 0.007 \ , 0.415 \pm 0.007)$ \\
whether graph $\chi$ is 0 & $(0.502 \pm 0.022 \ , 0.603 \pm 0.054 \ , 0.078 \pm 0.015)$ \\
\end{tabular}
\caption{{\sf {\small
Summary of the accuracy measures (goodness of fitness) for the various quantities as machine-learnt by a decision-tree classifier or simple feed-forward neural network within a 5-fold cross-validation.
The triple applies to cases of binary classification and refers to (naive precision, F1-score, Matthews $\phi$ coefficient) and the pair applies to multi-category classification and refers to (naive precision, Matthews $\phi$ coefficient).
}}
\label{t:summary}}
\end{table}

Thus prepared, we studied more precise bounds on the (random-walk normalized) Laplacian \cite{ly,ly2} which are analogues of the classical Li-Yau \cite{LiYau} inequalities for Riemannian manifolds.
We studied the statistics of the eigenvalue distributions in comparison with the various inequalities and then applied principal component analysis (PCA) and topological data analysis (TDA) on this distribution.
We found that planar and non-planar graphs have rather distinct principal components in the Euclidean point cloud of ordered Laplacian eigenvalues.
Then, using supervised learning on the adjacency matrices, a fairly satisfactory prediction of the spectral gap using the neural network was achieved.

Continuing in this vain, bearing in mind the intimate relation between curvature, Laplacian zero-modes and (co-)homology for manifolds, we proceeded to ask the question whether ML could distinguish Ricci-flat graphs.
There are two notions of Ricci-flatness for (locally) finite graphs, in terms of OLLY and of curvature-dimension inequalities.
We found that supervised ML can achieve more than 90\% accuracy (see Table \ref{t:summary} for the precise measures of goodness of fit) for both notions in a 5-fold cross-validation.
This is very assuring indeed.
It should be emphasized that all ML calculations involved in this paper are performed on a standard laptop using Mathematica, in a matter of seconds or minutes.
Finally, we explored the homology of graphs in the sense of \cite{glmy1,glmy2}.
We presented a rough distribution of Euler number over our database and attempted to machine-learn whether a graph is Euler number 0, though here the accuracies are not high and the ML did not perform much better than random guessing.

We hope our preliminary investigations have paved the road for countless new ventures in applying machine-learning and data scientific techniques to graphs and manifolds.
For example, we have not addressed the myriad of key functions in graph theory such as chromatic polynomials, Ihara zeta functions, etc.
In addition to the connections between graphs and manifolds, the study of digraphs in the guise of quivers is a hotly pursued topic in mathematics and physics ranging from cluster algebras to supersymmetric quantum field theories, to string theory (cf.~\cite{He:2011ge} for graph zeta function and gauge theory as well as machine-learning quiver gauge theories and cluster mutation \cite{Bao:2020nbi}).
Ultimately, we would like to see where in the hierarchy of complexity do combinatorial theorems on graphs reside as far as AI/ML is concerned and indeed whether new conjectures can be formulated \cite{hetalk}.
The terra incognita where ML meets pure mathematics beckons her ever-alluring invitation.

\section*{Acknowledgments}
We are grateful to Shuliang Bai and Yong Lin for their careful reading of and valuable comments on preliminary drafts of the paper.
Y.-H.~H is indebted to the Science and Technology Facilities Council, UK, for grant ST/J00037X/1 as well as a chair professorship from Nankai University where this work began.
The work of S.-T. Y. is supported in part by a grant from the Simons Foundation in Homological Mirror Symmetry.

\appendix
\section{Illustrative Example}
In this appendix, let us take a specific example from our database, and compute all the relevant quantities discussed throughout the main body.
The particulars thus detailed should serve an illustrative purpose.

Take the so-called {\bf dipyramid graph} $G_{dp}(5)$ on 5 vertices, whose figure, adjacency and degree matrices are shown as follows:
\begin{equation}
\begin{array}{c}
\includegraphics[trim=0mm 0mm 0mm 0mm, clip, width=2in]{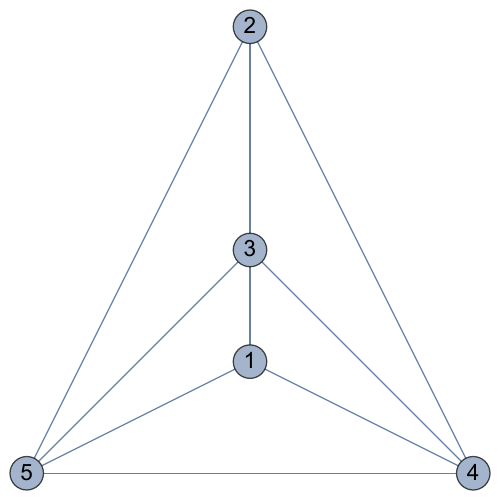}
\end{array}
A = {\scriptsize \left(
\begin{array}{ccccc}
 0 & 0 & 1 & 1 & 1 \\
 0 & 0 & 1 & 1 & 1 \\
 1 & 1 & 0 & 1 & 1 \\
 1 & 1 & 1 & 0 & 1 \\
 1 & 1 & 1 & 1 & 0 \\
\end{array}
\right)} \ , \quad
D = \mbox{Diag}(3, 3, 4, 4, 4) \ ,
\end{equation}
where the vertices have been labeled explicitly.
The Laplacian $L$ and the random-walk normalized Laplacian $\Delta$ are subsequently
\begin{equation}
L = D - A = 
{\scriptsize
\left(
\begin{array}{ccccc}
 3 & 0 & -1 & -1 & -1 \\
 0 & 3 & -1 & -1 & -1 \\
 -1 & -1 & 4 & -1 & -1 \\
 -1 & -1 & -1 & 4 & -1 \\
 -1 & -1 & -1 & -1 & 4 \\
\end{array}
\right)} \ , 
\quad
\Delta = D^{-1}L = 
{\scriptsize
\left(
\begin{array}{ccccc}
 1 & 0 & -\frac{1}{3} & -\frac{1}{3} & -\frac{1}{3} \\
 0 & 1 & -\frac{1}{3} & -\frac{1}{3} & -\frac{1}{3} \\
 -\frac{1}{4} & -\frac{1}{4} & 1 & -\frac{1}{4} & -\frac{1}{4} \\
 -\frac{1}{4} & -\frac{1}{4} & -\frac{1}{4} & 1 & -\frac{1}{4} \\
 -\frac{1}{4} & -\frac{1}{4} & -\frac{1}{4} & -\frac{1}{4} & 1 \\
\end{array}
\right)
} \ .
\end{equation}
The eigenvalues of $\Delta$ are, therefore, $\left\{\frac{3}{2},\frac{5}{4},\frac{5}{4},1,0\right\}$.

Clearly, the skewness of $G_{dp}(5)$ is 0 and the graph is planar; so too is the genus 0.
The chromatic number is 4 and it requires 4 colours, say with vertices $1,3,4,5$ with 4 different colours and vertex $2$ the same as 1.
The diameter is 2, exemplified by the minimal distance from vertex 1 to 2.
The girth is 3, since the minimal cycle is length 3.
There are no Euler cycles but there are 6 Hamiltonian cycles: $\big( \{1, 3, 2, 4, 5\}, \{1, 3, 2, 5, 4\}, \{1, 3, 4, 2, 5\}, \{1, 3, 5, 2, 4\}, \{1, 4, 2, 3, 5\}, \{1, 4, 3, 2, 5\} \big)$.

We now assign a random direction to $G_{dp}(5)$ as follows, and for reference, we include all the elementary paths of each length:
\begin{equation}
\begin{array}{c}
\includegraphics[trim=0mm 0mm 0mm 0mm, clip, width=2in]{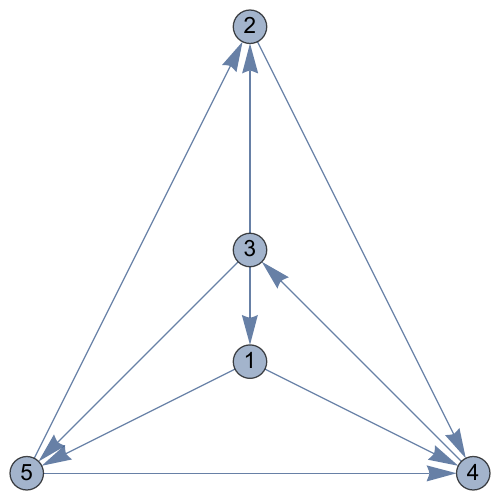}
\end{array}
\begin{array}{c|l}
p=1 & \{1,4\},\{1,5\},\{2,4\},\{3,1\},\{3,2\},\{3,5\},\{4,3\},\{5,2\},\{5,4\},\\
p=2 & \{1,4,3\},\{1,5,2\},\{1,5,4\},\{2,4,3\},\{3,1,4\},\{3,1,5\},\{3,2,4\}, \\
	& \{3,5,2\},\{3,5,4\},\{4,3,1\},\{4,3,2\},\{4,3,5\},\{5,2,4\},\{5,4,3\},\\
p=3 & \{1,4,3,2\},\{1,4,3,5\},\{1,5,2,4\},\{1,5,4,3\},\{2,4,3,1\},\\
	& \{2,4 ,3,5\},\{3,1,5,2\},\{3,1,5,4\},\{3,5,2,4\},\{4,3,1,5\},\\
	& \{4,3,5,2\},\{5,2,4,3\},\{5,4,3,1\},\{5,4,3,2\},\\
p=4 & \{1,4,3,5,2\},\{1,5,2,4,3\},\{1,5,4,3,2\},\{2,4,3,1,5\},\\
	& \{3,1,5,2,4\},\{4,3,1,5,2\},\{5,2,4,3,1\}
\end{array}
\end{equation}
Now, whilst \cite{glmy1,glmy2} give certain theorems on computing the (co-)homologies for classes of graphs, let us approach the problem by brute force for simplicity.

Let us consider $\cA_2$, whose basis is given (we introduce formal variables $x$ for convenience) as
\begin{align}
\nn
\cA_2 = \mbox{Span} & \left\{x_{\{1,4,3\}},x_{\{1,5,2\}},x_{\{1,5,4\}},x_{\{2,4,3
   \}},x_{\{3,1,4\}},x_{\{3,1,5\}},x_{\{3,2,4\}},x_{\{3,5,2
   \}},x_{\{3,5,4\}},x_{\{4,3,1\}}, \right. \\
   &\left. x_{\{4,3,2\}},x_{\{4,3,5
   \}},x_{\{5,2,4\}},x_{\{5,4,3\}}\right\}
\end{align}
from which we can compute the action of $\partial$ on each of the regular elementary 2-paths as:
\begin{align}
\nn
\partial \cA_2 = \mbox{Span} & \big\{
 -x_{\{1,3\}}+x_{\{1,4\}}+x_{\{4,3\}} \ ,
 -x_{\{1,2\}}+x_{\{1,5\}}+x_{\{5,2\}} \ ,
 -x_{\{1,4\}}+x_{\{1,5\}}+x_{\{5,4\}} \ ,
\\ \nn
& -x_{\{2,3\}}+x_{\{2,4\}}+x_{\{4,3\}} \ ,
 x_{\{1,4\}}+x_{\{3,1\}}-x_{\{3,4\}} \ ,
 x_{\{1,5\}}+x_{\{3,1\}}-x_{\{3,5\}} \ ,
 \\ \nn
& x_{\{2,4\}}+x_{\{3,2\}}-x_{\{3,4\}} \ ,
 -x_{\{3,2\}}+x_{\{3,5\}}+x_{\{5,2\}} \ ,
 -x_{\{3,4\}}+x_{\{3,5\}}+x_{\{5,4\}} \ ,
\\ \nn
& x_{\{3,1\}}-x_{\{4,1\}}+x_{\{4,3\}} \ , 
 x_{\{3,2\}}-x_{\{4,2\}}+x_{\{4,3\}} \ ,
 x_{\{3,5\}}+x_{\{4,3\}}-x_{\{4,5\}} \ ,
 \\ \nn
& x_{\{2,4\}}+x_{\{5,2\}}-x_{\{5,4\}} \ ,
 x_{\{4,3\}}-x_{\{5,3\}}+x_{\{5,4\}}
\big\} \ .
\end{align}
These need to be compared to the basis for $\cA_1$, which are
\begin{equation}
\cA_1 = \mbox{Span} 
\left\{x_{\{1,4\}},x_{\{1,5\}},x_{\{2,4\}},x_{\{3,1\}},x_{\{3,2\}},x_{\{3,5
   \}},x_{\{4,3\}},x_{\{5,2\}},x_{\{5,4\}}\right\} \ .
\end{equation}

We see that many variables appear in $\partial \cA_2$ which are {\em not} in $\cA_1$, which means that as it stands, the basis for $\partial \cA_2$ is not $\partial$-invariant.
In particular, these ``bad'' variables are
\begin{equation}
x_{bad} =
\left\{x_{\{1,2\}},x_{\{1,3\}},x_{\{2,3\}},x_{\{3,4\}},x_{\{4,1\}},x_{\{4,2
   \}},x_{\{4,5\}},x_{\{5,3\}}\right\} \ .
\end{equation}
The question then is: what linear combinations of $\partial \cA_2$, if any, exists so that when expanded out, all coefficients of the bad coefficients vanish? The dimenion of such linear combination of solutions is then the dimension of $\Omega_2$, the requisite $\partial$-invariant subspace of $\cA_2$.
There are 14 terms in $\partial \cA_2$ so we can set 14 arbitrary coefficients, imposing that when expanded the coefficients of the 8 $x_{bad}$ vanish gives a linear system which solves explicitly to give us 6 free coefficients, which means that $\dim \Omega_2 = 6$.
We have encountered the generic situation by coincidence, in general, we have to solve the linear system case by case.

We now repeat the above for $\cA_3, \cA_4$, etc.
In summary, we find that, in addition to the two by construction, viz.,
$\dim \Omega_0 = \dim \cA_0 = |V| = 5$ and
$\dim \Omega_1 = \dim \cA_1 = |E| = 9$, 
\begin{equation}
\dim \Omega_2 = 6, \ \dim \Omega_3 = 2, \ \dim \Omega_{n\geq 4} = 0, \ 
\end{equation}
so that the Euler number is
\begin{equation}
\chi(G_{dp}(5)) = 5 - 9 + 6 - 2 = 0 \ .
\end{equation}

\newpage 

\end{document}